\documentclass[12pt,twoside]{article}
\usepackage[hmargin=0.8in,vmargin=0.9in]{geometry}
\usepackage{fancyhdr}
\usepackage{graphicx}
\usepackage{subfigure}
\usepackage{amssymb}
\usepackage{amsmath}
\usepackage{amsfonts}
\usepackage{theorem}
\usepackage{mathrsfs}
\usepackage{mathtools}
\usepackage{bm}
\usepackage{color}
\usepackage{setspace}
\usepackage{exscale}
\usepackage{relsize}
\usepackage{epstopdf}
\usepackage{float}
\usepackage{cite}
\usepackage{nicefrac}
\usepackage{setspace}

\usepackage{booktabs,multirow} 
\usepackage{array} 
\usepackage{paralist} 
\usepackage{verbatim} 
\usepackage{subfigure} 

\theoremstyle{plain}			
\newtheorem{thm}{Theorem}[section]

\newtheorem{rmk}[thm]{Remark}

{\theorembodyfont{\rmfamily}}

\usepackage{tikz}
\usetikzlibrary{positioning}
\usepackage{xcolor}

\allowdisplaybreaks[1]

\numberwithin{equation}{section}
\numberwithin{figure}{section}
\numberwithin{table}{section}

\newcommand\eref[1]{(\ref{#1})}

\newcommand*\xbar[1]{%
  \hbox{%
    \vbox{%
      \hrule height 0.5pt 
      \kern0.4ex
      \hbox{%
        \kern-0.05em
        \ensuremath{#1}%
        \kern-0.00em
      }%
    }%
  }%
}

\newcommand{\bmF}{\bm{\mathcal{F}}}
\newcommand{\bmG}{\bm{\mathcal{G}}}
\newcommand{\mF}{\bm{F}}

\newcommand{\mG}{\bm{G}}

\newcommand{\mU}{\bm{U}}

\newcommand{\dt}{\Delta t}
\newcommand{\dx}{\Delta x}
\newcommand{\dy}{\Delta y}

\newcommand{\hf}{{\frac{1}{2}}}

\newcommand{\jph}{{j+\frac{1}{2}}}
\newcommand{\jmh}{{j-\frac{1}{2}}}
\newcommand{\kph}{{k+\frac{1}{2}}}
\newcommand{\kmh}{{k-\frac{1}{2}}}

\newcommand{\ajphp}{{a_{j+\frac{1}{2}}^+}}
\newcommand{\ajphm}{{a_{j+\frac{1}{2}}^-}}

\def\softd{{\leavevmode\setbox1=\hbox{d}%
          \hbox to 1.05\wd1{d\kern-0.4ex{\char039}\hss}}}

\title{New Adaptive Low-Dissipation Central-Upwind Schemes}
\author{Shaoshuai Chu\thanks{Department of Mathematics and Shenzhen International Center for Mathematics, Southern University of Science and
Technology, Shenzhen, 518055, China; {\tt chuss2019@mail.sustech.edu.cn}}~~and Alexander Kurganov\thanks{Department of Mathematics, Shenzhen
International Center for Mathematics and Guangdong Provincial Key Laboratory of Computational Science and Material Design, Southern
University of Science and Technology, Shenzhen, 518055, China; {\tt alexander@sustech.edu.cn}}}

\begin{document}

\date{}
\maketitle

\begin{abstract}
We introduce new second-order adaptive low-dissipation central-upwind (LDCU) schemes for the one- and two-dimensional hyperbolic systems of
conservation laws. The new adaptive LDCU schemes employ the LDCU numerical fluxes (recently proposed in [{\sc A. Kurganov and R. Xin},
J. Sci. Comput., 96 (2023), Paper No. 56]) computed using the point values reconstructed with the help of adaptively selected nonlinear
limiters. To this end, we use a smoothness indicator to detect ``rough'' parts of the computed solution, where the piecewise linear
reconstruction is performed using an overcompressive limiter, which leads to extremely sharp resolution of shock and contact waves. In the
``smooth'' areas, we use a more dissipative limiter to prevent appearance of artificial kinks and staircase-like structures there. In order
to avoid oscillations, we perform the reconstruction in the local characteristic variables obtained using the local characteristic
decomposition. We test two different smoothness indicators and apply the developed schemes to the one- and two-dimensional Euler equations
of gas dynamics. The obtained numerical results clearly demonstrate that the new adaptive LDCU schemes outperform the original ones.
\end{abstract}



\section{Introduction}
This paper focuses on developing new adaptive numerical methods for the hyperbolic systems of conservation laws, which in the one- (1-D) and
two-dimensional (2-D) cases, read as
\begin{equation}
\mU_t+\mF(\mU)_x=\bm0,
\label{1.1}
\end{equation}
and
\begin{equation}
\mU_t+\mF(\mU)_x+\mG(\mU)_y=\bm0,
\label{1.2}
\end{equation}
respectively. Here, $x$ and $y$ are spatial variables, $t$ is the time, $\mU\in\mathbb R^d$ is a vector of unknown functions, and
$\mF:\mathbb R^d\to\mathbb R^d$ and $\mG:\mathbb R^d\to\mathbb R^d$ are nonlinear fluxes.

It is well-known that even when the initial data are smooth, solutions of \eref{1.1} and \eref{1.2} can produce extremely complex nonsmooth
wave patterns including shocks, rarefactions, and contact discontinuities. This makes it quite challenging to develop accurate and reliable
shock-capturing numerical methods for \eref{1.1} and \eref{1.2}.

A library of numerical methods for the studied systems have been introduced since the pioneering works of Friedrichs \cite{Fri}, Lax
\cite{Lax}, and Godunov \cite{Godunov59}. We refer the reader to the monographs and review papers
\cite{KLR20,Tor,Leveque02,Hesthaven18,BAF,Shu20} and references therein, where one can find a description of many existing numerical
methods. In this paper, we restrict our consideration to semi-discrete finite-volume (FV) methods, where the solution, represented in terms
of its cell averages, is evolved in time with the help of the numerical fluxes, computed, in turn, using the reconstructed point values of
$\mU$ at the boundaries of the FV cells. Many of such schemes are upwind in the sense that their numerical fluxes are based on either exact
or approximate solution of the (generalized) Riemann problems arising at each cell interface. We, however, focus on the Riemann-solver-free
central-upwind (CU) schemes, which provide one with accurate, efficient and robust tools for a wide variety of hyperbolic systems. The CU
schemes belong to the class of non-oscillatory central schemes, but they have a certain upwind nature as they rely on the local one-sided
speeds of propagation, which can be estimated using the largest and smallest eigenvalues of the corresponding Jacobians. The original CU
schemes from \cite{Kurganov01,Kurganov02} contain relatively large amount of numerical dissipation, which was reduced in \cite{Kurganov07}
and recently in \cite{KX_22}, where built-in ``anti-diffusion'' terms were introduced. The amount of numerical dissipation can be also
reduced by applying the local characteristic decomposition (LCD) technique to the numerical diffusion of the CU fluxes; see \cite{CCHKL_22}.

In this paper, we use the low-dissipation CU (LDCU) numerical fluxes from \cite{KX_22}, and further enhance the resolution of the ``rough''
parts of the computed solution by applying a new scheme adaption approach: The point values used to evaluate the LDCU fluxes are
reconstructed with the help of adaptively selected nonlinear limiters, which are, in general, required to make the reconstructed point
values non-oscillatory. A variety of limiters are available; see, e.g., \cite{Lie03,Nessyahu90,Sweby84,BAF,Hesthaven18,Leveque02,Tor} and
references therein. Many of the limiters can be classified as dissipative, compressive, or overcompressive as it was done in \cite{Lie03}.
The use of compressive and overcompressive limiters leads to very sharp resolution of discontinuous parts of the approximated solution,
while dissipative limiters may smear the jumps. At the same time, applying compressive and overcompressive limiters in the smooth areas
typically results in the artificial sharpening of the smooth solution profiles, that is, in the appearance of kinks or staircase-like
structures, or even non-physical jumps.

We therefore switch between different limiters. To this end, we need to automatically  detect ``rough'' (nonsmooth) parts of the computed
solution with the help of a smoothness indicator (SI). Many different SIs are readily available; see, e.g.,
\cite{DZLD14,Dewar15,ABD08,GT06,GT02,GPP,PupSem,QS05,VR16,FS17} and references therein. In this paper, we test two different SIs: a slightly
modified minmod (MM)-based shock indicator from \cite{WSYK15} (see also \cite{Harten89,SO89}) and a SI based on weak local residuals (WLR)
from \cite{Kurganov12a} (see also \cite{KK05,Dewar15,KKP02}). In the areas identified as being ``rough'', we use the overcompressive SBM
limiters from \cite{Lie03}, while switching to the dissipative Minmod2 limiter elsewhere. It is well-known that the use of any of these two
limiters may lead to numerical oscillations in the vicinities of shock and contact discontinuities. In order to reduce these oscillations,
we perform the reconstruction in the local characteristic variables rather than in the conservative or primitive ones (this strategy was
advocated in, e.g., \cite{Qiu02}). We switch to the characteristic variables using the LCD, which is often used in the context of high-order
schemes, but can also be implemented to enhance the resolution of second-order schemes; see, e.g., \cite{CCHKL_22,Joh,Qiu02,Shu20} and
references therein.

The paper is organized as follows. In \S\ref{sec2}, we review the recently proposed 1-D LDCU scheme from \cite{KX_22}. We then introduce the
adaptive schemes that employ either the MM- or WLR-based SI to detect the ``rough'' areas. In \S\ref{sec3}, we extend the proposed adaptive
LDCU scheme to the 2-D case. In \S\ref{sec4}, we apply the developed schemes to a number of 1-D and 2-D numerical examples for the Euler
equations of gas dynamics. We demonstrate that the adaptive LDCU schemes contain substantially smaller amount of numerical dissipation and
achieve much higher resolution compared with the LDCU schemes based on the Minmod2 limiters applied throughout the entire computational
domain. Finally, we give some concluding remarks in \S\ref{sec5}.

\section{One-Dimensional Scheme Adaption Algorithm}\label{sec2}
In this section, we consider the 1-D conservation laws \eref{1.1} and describe the 1-D adaptive algorithm.

\subsection{1-D Low-Dissipation Central-Upwind (LDCU) Schemes}
Assume that the computational domain is covered with the uniform cells $C_j:=[x_\jmh,x_\jph]$ with $x_\jph-x_\jmh\equiv\dx$ centered at
$x_j=(x_\jmh+x_\jph)/2$ and denote by $\xbar\mU_j(t)$ cell averages of $\mU(\cdot,t)$ over the corresponding intervals $C_j$, that is,
\begin{equation*}
\xbar\mU_j(t):\approx\frac{1}{\dx}\int\limits_{C_j}\mU(x,t)\,{\rm d}x.
\end{equation*}
We suppose that at a certain time $t\ge0$, the point values of the computed solution $\xbar \mU_j(t)$ are available. Note that all of the
indexed quantities are time-dependent, but from here on, we will suppress the time-dependence of all of the indexed quantities for the sake
of brevity.

According to the semi-discrete LDCU scheme from \cite{KX_22}, the computed cell averages are evolved in time by numerically solving the
following system of ordinary differential equations (ODEs):
\begin{equation}
\frac{{\rm d}\xbar\mU_j}{{\rm d}t}=-\frac{\bm{{\cal F}}_\jph-\bm{{\cal F}}_\jmh}{\dx},
\label{2.1}
\end{equation}
where $\bm{{\cal F}}_\jph$ are the LDCU numerical fluxes defined by
\begin{equation*}
\bm{{\cal F}}_\jph\big(\bm U_\jph^-,\bm U_\jph^+\big)=\frac{\ajphp\mF^-_\jph-\ajphm \mF^+_\jph}{\ajphp-\ajphm}+
\frac{\ajphp\ajphm}{\ajphp-\ajphm}\left(\mU^+_\jph-\mU^-_\jph\right)+\bm q_\jph.
\end{equation*}
Here, $\mF^\pm_\jph:=\mF\big(\mU^\pm_\jph\big)$ and $\mU^\pm_\jph$ are the right/left-sided point values of $\mU$ at the cell interface
$x=x_\jph$. The point values $\mU^\pm_\jph$ are reconstructed out of the given set of cell averages $\{\xbar\mU_j\}$ using a proper
nonlinear limiter; see \S\ref{sec211}. The one-sided local speeds  of propagation $a^\pm_\jph$ are estimated using the largest and the
smallest eigenvalues of the Jacobian $A(\mU):=\frac{\partial\mF}{\partial\mU}(\mU)$, $\lambda_1(A(\mU))\le\ldots\le\lambda_d(A(\mU))$. This
can be done, for example, by taking
\begin{equation*}
\begin{aligned}
&a^+_\jph=\max\big\{\lambda_d\big(A(\mU^+_\jph)\big),\lambda_d\big(A(\mU^-_\jph)\big),0\big\},\\
&a^-_\jph=\min\big\{\lambda_1\big(A(\mU^+_\jph)\big),\lambda_1\big(A(\mU^-_\jph)\big),0\big\}.
\end{aligned}
\end{equation*}

Finally, $\bm q_\jph$ is a built-in ``anti-diffusion'' term, which can be derived for a particular system at hand. For instance, we consider
the 1-D Euler equations of gas dynamics, which read as \eref{1.1} with
\begin{equation}
\mU=\big(\rho,\rho u,E\big)^\top\quad{\rm and}\quad\mF=\big(\rho u,\rho u^2+p,u(E+p)\big)^\top.
\label{2.4}
\end{equation}
Here, $\rho$, $u$, $p$, and $E$ are the density, velocity, pressure, and total energy, respectively, and the system is completed through the
following equations of state (EOS) for ideal gases:
\begin{equation}
p=(\gamma-1)\Big[E-\hf\rho u^2\Big],
\label{2.9}
\end{equation}
where the parameter $\gamma$ represents the specific heat ratio. For the Euler system \eref{1.1}, \eref{2.4}, \eref{2.9}, the
``anti-diffusion'' term $\bm q_\jph$ has been rigorously derived in \cite{KX_22} and it is given by
\begin{equation*}
\bm q_\jph={\rm minmod}\big(-a^-_\jph(\rho^*_\jph-\rho^-_\jph),a^+_\jph(\rho^+_\jph-\rho^*_\jph)\big)
\begin{pmatrix}1\\u^*_\jph\\[0.9ex]\hf\big(u^*_\jph\big)^2\end{pmatrix}.
\end{equation*}
Here, $\rho^*_\jph$ and $(\rho u)^*_\jph$ are the first and second components of
\begin{equation*}
\mU^*_\jph=\frac{a^+_\jph\mU^+_\jph-a^-_\jph\mU^-_\jph-\left\{\bm F(\mU^+_\jph)-\bm F(\mU^-_\jph)\right\}}{a^+_\jph-a^-_\jph},
\end{equation*}
$u^*_\jph=(\rho u)^*_\jph/\rho^*_\jph$, and the minmod function is defined by
\begin{equation*}
{\rm minmod}(z_1,z_2,\ldots):=\begin{cases}
\min_j\{z_j\}&\mbox{if}~z_j>0\,\,\forall\,j,\\
\max_j\{z_j\}&\mbox{if}~z_j<0\,\,\forall\,j,\\
0            &\text{otherwise.}
\end{cases}
\end{equation*}

\subsubsection{Nonlinear Limiters}\label{sec211}
As mentioned before, the point values $\mU^\pm_\jph$ are obtained with the help of a conservative piecewise linear reconstruction, designed
using a proper nonlinear limiter. In this paper, we use a family of the SBM limiters (introduced in \cite{Lie03}) applied to the local
characteristic variables. To this end, we introduce the matrices $\widehat A_\jph=A\big((\,\xbar \mU_j+\xbar \mU_{j+1})/2\big)$ and compute
the matrices $R_\jph$ and $R^{-1}_\jph$ such that $R^{-1}_\jph\widehat A_\jph R_\jph$ are diagonal matrices. We then introduce the local
characteristic variables $\bm\Gamma$ in the neighborhood of $x=x_\jph$:
$$
\bm\Gamma_k=R^{-1}_\jph\xbar\mU_k,\quad k=j-1,\,j,\,j+1,\,j+2.
$$
Equipped with the values $\bm\Gamma_{j-1}$, $\bm\Gamma_j$, $\bm\Gamma_{j+1}$, and $\bm\Gamma_{j+2}$, we compute the slopes
\begin{equation}
(\bm\Gamma_x)_j=\phi^{\rm SBM}_{\theta,\tau}\left(\frac{\bm\Gamma_{j+1}-\bm\Gamma_j}{\bm\Gamma_j-\bm\Gamma_{j-1}}\right)
\frac{\bm\Gamma_j-\bm\Gamma_{j-1}}{\dx}
\label{2.10}
\end{equation}
and
\begin{equation}
(\bm\Gamma_x)_{j+1}=\phi^{\rm SBM}_{\theta,\tau}\left(\frac{\bm\Gamma_{j+2}-\bm\Gamma_{j+1}}{\bm\Gamma_{j+1}-\bm\Gamma_j}\right)
\frac{\bm\Gamma_{j+1}-\bm\Gamma_j}{\dx},
\label{2.11}
\end{equation}
where the two-parameter SBM function
\begin{equation}
\phi^{\rm SBM}_{\theta,\tau}(r):=\begin{cases}0&\mbox{if $r<0$},\\\min\{r\theta,1+\tau(r-1)\}&\mbox{if}~0<r\le1,\\
r\phi^{\rm SBM}_{\theta,\tau}(\frac{1}{r})&\text{otherwise,}\end{cases}
\label{2.12}
\end{equation}
is applied in the component-wise manner.

The parameters $\theta\in[1,2]$ and $\tau$ in \eref{2.12} can be used to control the amount of numerical dissipation present in the
resulting scheme. First, larger $\theta$'s correspond to less dissipative but, in general, more oscillatory reconstructions. In all of the
numerical examples reported in \S\ref{sec4}, we have taken $\theta=2$. Second, according to \cite{Lie03}, if $\tau\ge0.5$, then the SBM
limiter is dissipative and its use typically causes contact discontinuities to be severely smeared in time. If $0\le\tau<0.5$, then the SBM
limiter is compressive and in this case, contact waves are usually resolved sharply within few points, but smooth extrema might be slightly
compressed resulting in continuous solution profiles having a kink. If $\tau<0$, then the limiter is overcompressive so that contact
discontinuities typically stay very sharp for long time, while smooth solutions become overcompressed as time evolves resulting in the
appearance of artificial ${\cal O}(1)$ jump discontinuities.

Equipped with \eref{2.10} and \eref{2.11}, we evaluate
$$
\bm\Gamma^-_\jph=\bm\Gamma_j+\frac{\dx}{2}(\bm\Gamma_x)_j\quad\mbox{and}\quad
\bm\Gamma^+_\jph=\bm\Gamma_{j+1}-\frac{\dx}{2}(\bm\Gamma_x)_{j+1},
$$
and then obtain the corresponding point values of $\mU$ by
\begin{equation*}
\mU^\pm_\jph=R_\jph\bm\Gamma^\pm_\jph.
\end{equation*}
\begin{rmk}
For detailed explanations on how the matrices $R_\jph$ and $R^{-1}_\jph$ are computed in the case of the Euler equation of gas dynamics, we
refer the reader to \cite[Appendix A]{CCHKL_22}.
\end{rmk}

\subsection{One-Dimensional Adaptive Schemes}
We now turn to the description of the proposed adaptive schemes. The key ingredient of the new schemes is the use of the different limiters
from the family \eref{2.12} in different parts of the computational domain. In particular, we use an overcompressive limiter with
$\tau=-0.25$ in the ``rough'' parts of the computed solution and a dissipative limiter with $\tau=0.5$ elsewhere. The latter limiter is, in
fact, the Minmod2 limiter, which can be written in a simpler form since \eref{2.10} and \eref{2.11} with $\theta=2$ and $\tau=0.5$ reduce to
\begin{equation*}
(\bm\Gamma_x)_j={\rm minmod}\left(2\,\frac{\bm\Gamma_j-\bm\Gamma_{j-1}}{\dx},\,\frac{\bm\Gamma_{j+1}-\bm\Gamma_{j-1}}{2\dx},\,
2\,\frac{\bm\Gamma_{j+1}-\bm\Gamma_j}{\dx}\right),
\end{equation*}
and
\begin{equation*}
(\bm\Gamma_x)_{j+1}={\rm minmod}\left(2\,\frac{\bm\Gamma_{j+1}-\bm\Gamma_j}{\dx},\,\frac{\bm\Gamma_{j+2}-\bm\Gamma_j}{2\dx},\,
2\,\frac{\bm\Gamma_{j+2}-\bm\Gamma_{j+1}}{\dx}\right).
\end{equation*}

In order to implement this simple scheme adaption approach, we need to automatically detect ``rough'' parts of the computed solution. This
is done using either the MM- or WLR-based SIs briefly described in \S\ref{sec2.2.1} and \S\ref{sec2.2.2} below.

\subsubsection{Minmod-Based Smoothness Indicator}\label{sec2.2.1}
We first compute the MM-based quantities
$$
s_j={\rm minmod}\big(\,\xbar\rho_{j+1}-\xbar\rho_j,\,\xbar\rho_j-\xbar\rho_{j-1}\big),
$$
and then we say that the cell $C_j$ is ``rough'' if $|s_j|>\max\big\{|s_{j-1}|,|s_{j+1}|\big\}+\delta$, where $\delta$ is a small positive
number. The parameter $\delta$ has to be selected for each problem at hand and it should indicate a size of a jump in $\rho$, which we
neglect when detecting ``rough'' parts of the solution. In fact, this SI is not very sensitive to the choice of $\delta$ and in all of the
numerical examples reported in \S\ref{sec4}, we have taken $\delta=10^{-4}$.

\subsubsection{Weak Local Residual-Based Smoothness Indicator}\label{sec2.2.2}
In order to detect ``rough'' areas, one can also use the WLR-based SI, which we obtain as follows. First, we assume that the cell averages
$\xbar\mU_j$ are available at a certain time level $t=t^n$ and the two previous time levels $t=t^{n-1}$ (with $t^n-t^{n-1}=\dt^{n-1}$) and
$t^{n-2}$ (with $t^{n-1}-t^{n-2}=\dt^{n-2}$). In addition, we assume that the solution has been reconstructed at $t=t^{n-1}$ and $t^{n-2}$
and the corresponding point values at the cell interfaces $x=x_\jph$ are available. We will denote these reconstructed point values by
$\mU^{n-1}_\jph$ and $\mU^{n-2}_\jph$. Recall that we, in fact, obtain two point values at each cell interface ($\mU^\pm_\jph$) and any of
them can be used to evaluate the WLRs for the density equation $\rho_t+(\rho u)_x=0$. We denote these WLRs by
$\varepsilon^{n-\frac{3}{2}}_\jph$ and compute them according to \cite{Kurganov12a}:
\begin{equation}
\begin{aligned}
\varepsilon^{n-\frac{3}{2}}_\jph&=\frac{\dx}{6}\left[\rho^{n-1}_{j+\frac{3}{2}}-\rho^{n-2}_{j+\frac{3}{2}}+
4\big(\rho^{n-1}_\jph-\rho^{n-2}_\jph\big)+\rho^{n-1}_\jmh-\rho^{n-2}_\jmh\right]\\
&+\frac{\dt^{n-2}}{4}\left[(\rho u)^{n-1}_{j+\frac{3}{2}}-(\rho u)^{n-1}_\jmh+(\rho u)^{n-2}_{j+\frac{3}{2}}-(\rho u)^{n-2}_\jmh\right].
\end{aligned}
\label{2.12a}
\end{equation}
The desired SIs are then obtained at each cell interface $x=x_\jph$ by setting
\begin{equation*}
\xbar\varepsilon^{\,n-\frac{3}{2}}_\jph:=\frac{1}{6}\left[\varepsilon^{\,n-\frac{3}{2}}_\jmh+4\varepsilon^{\,n-\frac{3}{2}}_\jph+
\varepsilon^{\,n-\frac{3}{2}}_{j+\frac{3}{2}}\right].
\end{equation*}
As mentioned in \cite{Kurganov12a}, the size of the WLRs and thus of the SIs for second-order schemes are expected to be
\begin{equation}
||\varepsilon^{n-\frac{3}{2}}||_\infty\sim
\begin{cases}\Delta&\mbox{near shock waves,}\\\Delta^\alpha&\mbox{near contact,}\\\Delta^4&\mbox{in smooth regions},\end{cases}
\label{2.13}
\end{equation}
where $\Delta:=\max\{\dt,\dx\}$ and $1<\alpha\le2$; see also \cite{KK05}.

Finally, we take advantage of \eref{2.13}, which suggests that the size of SIs ranges from ${\cal O}(\Delta)$ near shocks to
${\cal O}(\Delta^4)$ in the smooth regions and develop the following simple strategy for the automatic detection of ``rough'' parts of the
computed solution $\{\,\xbar\mU_j\}$. We mark the cell $C_j$ as ``rough'' as long as
\begin{equation}
\xbar \varepsilon^{\,n-\frac{3}{2}}_j>\texttt{C}(\dx)^2,
\label{2.14}
\end{equation}
where $\texttt{C}$ is a positive tunable constant to be selected for each problem at hand. The robustness of this shock detection strategy
depends on the sensitivity of the proposed algorithm to the choice of $\texttt{C}$. In principle, $\texttt{C}$ can be tuned on a coarse mesh
and then used for fine mesh computations, but as we demonstrate in Examples 2 in \S\ref{sec4}, this approach may fail. Therefore, even
though the use of the WLR-based SI may lead to extremely sharp results (like in, for instance, Example 3 in \S\ref{sec4}), the adaption
strategy that relies on this SI may not
be robust.
\begin{rmk}
While implementing \eref{2.12a}, we have used $\mU_\jph=\mU^-_\jph$ in all the numerical examples reported in \S\ref{sec4}. Other choices
like $\mU_\jph=\mU^+_\jph$ or $\mU_\jph=\big(\mU^-_\jph+\mU^+_\jph\big)/2$ can also be used with no visible advantages or disadvantages of
any of them.
\end{rmk}

\section{Two-Dimensional Scheme Adaption Algorithm}\label{sec3}
In this section, we extend the 1-D adaptive strategy introduced in \S\ref{sec2} to the 2-D hyperbolic systems of conservation laws
\eref{1.2}.

\subsection{2-D Low-Dissipation Central-Upwind (LDCU) Schemes}
Let the computational domain be covered with uniform cells $C_{j,k}:=[x_\jmh,x_\jph]\times[y_\kmh,y_\kph]$ with $x_\jph-x_\jmh\equiv\dx$ and
$y_\kph-y_\kmh\equiv\dy$ centered at $(x_j,y_k)$ with $x_j=(x_\jmh+x_\jph)/2$ and $y_k=(y_\kmh+y_\kph)/2$. We assume that the cell averages,
\begin{equation*}
\xbar\mU_{j,k}:\approx\frac{1}{\dx\dy}\iint\limits_{C_{j,k}}\mU(x,y,t)\,{\rm d}y\,{\rm d}x,
\end{equation*}
have been computed at a certain time $t\ge0$.

According to the semi-discrete LDCU scheme from \cite{KX_22}, the computed cell averages are evolved in time by numerically solving the
following system of ODEs:
\begin{equation}
\frac{{\rm d}\xbar\mU_{j,k}}{{\rm d}t}=-\frac{\bmF_{\jph,k}-\bmF_{\jmh,k}}{\dx}-\frac{\bmG_{j,\kph}-\bmG_{j,\kmh}}{\dy},
\label{3.1}
\end{equation}
where $\bm{{\cal F}}_{\jph,k}=\bm{{\cal F}}_{\jph,k}\big(\bm U_{\jph,k}^-,\bm U_{\jph,k}^+\big)$ and
$\bm{{\cal G}}_{j,\kph}=\bm{{\cal G}}_{j,\kph}\big(\bm U_{j,\kph}^-,\bm U_{j,\kph}^+\big)$ are the LDCU numerical fluxes defined by
\begin{equation*}
\begin{aligned}
\bmF_{\jph,k}&=\frac{a^+_{\jph,k}\mF^-_{\jph,k}-a^-_{\jph,k}\mF^+_{\jph,k}}
{a^+_{\jph,k}-a^-_{\jph,k}}+\frac{a^+_{\jph,k}a^-_{\jph,k}}{a^+_{\jph,k}-a^-_{\jph,k}}\left[\mU^+_{\jph,k}-\mU^-_{\jph,k}\right]+
\bm q^x_{\jph,k},\\
\bmG_{j,\kph}&=\frac{b^+_{j,\kph}\mG^-_{j,\kph}-b^-_{j,\kph}\mG^+_{j,\kph}}
{b^+_{j,\kph}-b^-_{j,\kph}}+\frac{b^+_{j,\kph}b^-_{j,\kph}}{b^+_{j,\kph}-b^-_{j,\kph}}\left[\mU^-_{j,\kph}-\mU^+_{j,\kph}\right]+
\bm q^y_{j,\kph}.
\end{aligned}
\end{equation*}
Here, $\mF^\pm_{\jph,k}:=\mF\big(\mU^\pm_{\jph,k}\big)$ and $\mG^\pm_{j,\kph}:=\mG\big(\mU^\pm_{j,\kph}\big)$, and $\mU^\pm_{\jph,k}$ and
$\mU^\pm_{j,\kph}$ are the one-sided point values of $\bm U$ at the cell interfaces $(x_\jph,y_k)$ and $(x_j,y_\kph)$, respectively. We
reconstruct the point values $\mU^\pm_{\jph,k}$ and $\mU^\pm_{j,\kph}$ using the LCD; see Appendix \ref{appa} for details. The one-sided
local speeds of propagation in the $x$- and $y$-directions, $a^\pm_{\jph,k}$ and $b^\pm_{j,\kph}$, can be estimated by the largest and
smallest eigenvalues of the Jacobians $A(\mU):=\frac{\partial\mF}{\partial\mU}(\mU)$ and $B(\mU):=\frac{\partial\mG}{\partial\mU}(\mU)$, for
example, by setting
\begin{equation*}
\begin{aligned}
&a^+_{\jph,k}=\max\left\{\lambda_d\big(A(\mU^+_{\jph,k})\big),\lambda_d\big(A(\mU^-_{\jph,k})\big),0\right\},\\
&a^-_{\jph,k}=\min\left\{\lambda_1\big(A(\mU^+_{\jph,k})\big),\lambda_1\big(A(\mU^-_{\jph,k})\big),0\right\},\\
&b^+_{j,\kph}=\max\left\{\lambda_d\big(B(\mU^+_{j,\kph})\big),\lambda_d\big(B(\mU^-_{j,\kph})\big),0\right\},\\
&b^-_{j,\kph}=\min\left\{\lambda_1\big(B(\mU^+_{j,\kph})\big),\lambda_1\big(B(\mU^-_{j,\kph})\big),0\right\}.
\end{aligned}
\end{equation*}

Finally, $\bm q^x_{\jph,k}$ and $\bm q^y_{j,\kph}$ are built-in ``anti-diffusion'' terms, which can be derived for a particular system
\eref{1.2} at hand. For instance, we consider the 2-D Euler equations of gas dynamics, which read as \eref{1.2} with
\begin{equation}
\mU=\big(\rho,\rho u,\rho v,E\big)^\top,~~\mF=\big(\rho u,\rho u^2+p,\rho uv,u(E+p)\big)^\top,~~
\mG=\big(\rho v,\rho uv,\rho v^2+p,v(E+p)\big)^\top,
\label{3.3}
\end{equation}
where $v$ is the $y$-velocity and the other variables are as the same as in \eref{2.4}. The system \eref{3.3} is completed through the
following EOS for ideal gases:
\begin{equation}
p=(\gamma-1)\Big[E-\frac{\rho}{2}(u^2+v^2)\Big].
\label{3.4}
\end{equation}
For the sake of brevity, we omit the details on the built-in ``anti-diffusion'' terms $\bm q^x_{\jph,k}$ and $\bm q^y_{j,\kph}$. For the
system \eref{1.2}, \eref{3.3}--\eref{3.4}, they have been derived in \cite{KX_22}.

\subsection{Two-Dimensional Adaptive Schemes}
We now turn to the description of the proposed adaptive schemes for the 2-D system. As in the 1-D case, we use an overcompressive SBM
limiter with $\tau=-0.25$ in the ``rough'' parts of the computed solution and a dissipative Minmod2 limiter elsewhere. To this end, we
detect the ``rough'' parts of the numerical solution using either the MM- or WLR-based SIs described in \S\ref{sec3.2.1} and
\S\ref{sec3.2.2} below.

\subsubsection{Two-Dimensional Minmod-Based Smoothness Indicator}\label{sec3.2.1}
The 1-D MM-based SI introduced in \S\ref{sec2.2.1}, is extended to the 2-D case in the ``dimension-by-dimension'' manner. We first compute
the MM-based quantities in the $x$-direction,
$$
s^x_{j,k}={\rm minmod}\big(\rho_{j+1,k}-\rho_{j,k},\rho_{j,k}-\rho_{j-1,k}\big),
$$
and use the overcompressive SBM limiter to compute the slopes in the $x$-direction only in those cell $C_{j,k}$, where
$|s^x_{j,k}|>\max\big\{|s^x_{j-1,k}|,|s^x_{j+1,k}|\big\}+\delta$. Similarly, we compute the MM-based quantities in the $y$-direction,
$$
s^y_{j,k}={\rm minmod}\big(\rho_{j,k+1}-\rho_{j,k},\rho_{j,k}-\rho_{j,k-1}\big),
$$
and use the overcompressive SBM limiter to compute the slopes in the $y$-direction only in those cell $C_{j,k}$, where
$|s^y_{j,k}|>\max\big\{|s^y_{j,k-1}|,|s^y_{j,k+1}| \big\}+\delta$.

\subsubsection{Two-Dimensional Weak Local Residual-Based Smoothness Indicator}\label{sec3.2.2}
One can also detect the ``rough'' parts of the numerical solution using the 2-D WLR-based SI, which we obtain as follows. First, we compute
the WLRs introduced in \cite{Kurganov12a}. For the 2-D density equation $\rho_t+(\rho u)_x+(\rho v)_y=0$, these WLRs are
\begin{equation*}
\varepsilon^{n-\frac{3}{2}}_{\jph,\kph}=\frac{1}{36\Delta}\dx\dy\,{\cal U}_{\jph,\kph}^{n-\frac{3}{2}}+
\frac{1}{12\Delta}\Big(\dy\dt{\cal F}_{j,k}^{n-\frac{3}{2}}+\dx\dt{\cal G}_{\jph,\kph}^{n-\frac{3}{2}}\Big),
\end{equation*}
where $\Delta:=\max\{\dt,\dx,\dy\}$ and
\begin{equation*}
\begin{aligned}
{\cal U}_{\jph,\,\kph}^{n-\frac{3}{2}}&=\left[\rho^{n-1}_{j+\frac{3}{2},k+\frac{3}{2}}-\rho^{n-2}_{j+\frac{3}{2},k+\frac{3}{2}}+
\rho^{n-1}_{j+\frac{3}{2},\kmh}-\rho^{n-2}_{j+\frac{3}{2},\kmh}+\rho^{n-1}_{\jmh,k+\frac{3}{2}}-\rho^{n-2}_{\jmh,k+\frac{3}{2}}\right.\\
&\left.+\rho^{n-1}_{\jmh,\kmh}-\rho^{n-2}_{\jmh,\kmh}\right]+4\left[\rho^{n-1}_{j+\frac{3}{2},\kph}-\rho^{n-2}_{j+\frac{3}{2},\kph}+
\rho^{n-1}_{\jmh,\kph}-\rho^{n-2}_{\jmh,\kph}\right.\\
&\left.+\rho^{n-1}_{\jph,k+\frac{3}{2}}-\rho^{n-2}_{\jph,k+\frac{3}{2}}+\rho^{n-1}_{\jph,\kmh}-\rho^{n-2}_{\jph,\kmh}\right]+
16\left[\rho^{n-1}_{\jph,\kph}-\rho^{n-2}_{\jph,\kph}\right],\\[0.8ex]
{\cal F}_{\jph,\kph}^{n-\frac{3}{2}}&=\left[(\rho u)^{n-1}_{j+\frac{3}{2},k+\frac{3}{2}}-(\rho u)^{n-1}_{\jmh,k+\frac{3}{2}}+
(\rho u)^{n-1}_{j+\frac{3}{2},\kmh}-(\rho u)^{n-1}_{\jmh,\kmh}\right.\\
&\left.+(\rho u)^{n-2}_{j+\frac{3}{2},k+\frac{3}{2}}-(\rho u)^{n-2}_{\jmh,k+\frac{3}{2}}+(\rho u)^{n-2}_{j+\frac{3}{2},\kmh}-
(\rho u)^{n-2}_{\jmh,\kmh}\right]\\
&+4\left[(\rho u)^{n-1}_{j+\frac{3}{2},\kph}-(\rho u)^{n-1}_{\jmh,\kph}-(\rho u)^{n-2}_{j+\frac{3}{2},\kph}-(\rho u)^{n-2}_{\jmh,\kph}
\right],\\[0.8ex]
{\cal G}_{\jph,\kph}^{n-\frac{3}{2}}&=\left[(\rho v)^{n-1}_{j+\frac{3}{2},k+\frac{3}{2}}-(\rho v)^{n-1}_{j+\frac{3}{2},\kmh}+
(\rho v)^{n-1}_{\jmh,k+\frac{3}{2}}-(\rho v)^{n-1}_{\jmh,\kmh}\right.\\
&\left.+(\rho v)^{n-2}_{j+\frac{3}{2},k+\frac{3}{2}}-(\rho v)^{n-2}_{j+\frac{3}{2},\kmh}+(\rho v)^{n-2}_{\jmh,k+\frac{3}{2}}-
(\rho v)^{n-2}_{\jmh,\kmh}\right]\\
&+4\left[(\rho v)^{n-1}_{\jph,k+\frac{3}{2}}-(\rho v)^{n-1}_{\jph,\kmh}-(\rho v)^{n-2}_{\jph,k+\frac{3}{2}}-(\rho v)^{n-2}_{\jph,\kmh}
\right].
\end{aligned}
\end{equation*}

Note that no values at the cell corners $(x_\jph,y_\kph)$ are reconstructed. We therefore set
$$
\mU^{\ell}_{\jph,\kph}:=\frac{1}{4}\left[\,\xbar\mU^{\,\ell}_{j,k}+\xbar\mU^{\,\ell}_{j+1,k}+\xbar\mU^{\,\ell}_{j,k+1}+
\xbar\mU^{\,\ell}_{j+1,k+1}\right],\quad\ell=n-2,\,n-1,
$$
where $\,\xbar\mU^{\,n-2}_{j,k}$ and $\,\xbar\mU^{\,n-1}_{j,k}$ denote the cell averages computed at times $t=t^{n-2}$ and $t=t^{n-1}$,
respectively.

The desired SIs are then obtained at each cell corner $(x_\jph,y_\kph)$ using a proper averaging, for instance, by
\begin{equation*}
\begin{aligned}
\xbar\varepsilon^{\,n-\frac{3}{2}}_{\jph,\kph}:&=\frac{1}{36}\left[\varepsilon^{n-\frac{3}{2}}_{\jmh,\kmh}+
\varepsilon^{n-\frac{3}{2}}_{\jmh,k+\frac{3}{2}}+\varepsilon^{n-\frac{3}{2}}_{j+\frac{3}{2},\kmh}+
\varepsilon^{n-\frac{3}{2}}_{j+\frac{3}{2},k+\frac{3}{2}}\right.\\
&\left.+\,4\left(\varepsilon^{n-\frac{3}{2}}_{\jmh,\kph}+\varepsilon^{n-\frac{3}{2}}_{\jph,\kmh}+
\varepsilon^{n-\frac{3}{2}}_{\jph,k+\frac{3}{2}}+\varepsilon^{n-\frac{3}{2}}_{j+\frac{3}{2},\kph}\right)+
16\,\varepsilon^{n-\frac{3}{2}}_{\jph,\kph}\right].
\end{aligned}
\end{equation*}

As mentioned in \cite{Kurganov12a}, the size of the WLRs and thus of the SIs for second-order schemes are expected to be the same as in
\eref{2.13}. We therefore act similarly to \eref{2.14} and mark the cell $C_{j,k}$ as ``rough'' as long as
$\,\xbar\varepsilon^{\,n-\frac{3}{2}}_{\jph,\kph}>\texttt{C}\max\{(\dx)^2,(\dy)^2\}$, where $\texttt{C}$ is a positive tunable constant to
be selected for each problem at hand.

\section{Numerical Examples}\label{sec4}
In this section, we test the developed adaptive schemes on several numerical examples. To this end, we compare the performance of the
original LDCU and the adaptive LDCU schemes by applying them to a number of initial-boundary value problems for the 1-D and 2-D Euler
equations of gas dynamics. The adaptive LDCU schemes with the MM- and WLR-based SIs used to detect the ``rough'' areas will be referred to
as the \emph{A-MM} and \emph{A-WLR} schemes, respectively.

In all of the numerical examples reported, we have solved the ODE systems \eref{2.1} and \eref{3.1} using the three-stage third-order strong
stability preserving (SSP) Runge-Kutta method; see, e.g., \cite{Gottlieb11,Gottlieb12}. We take $\gamma=1.4$ in Example 1--5 and
$\gamma=5/3$ in Example 6. The CFL number is 0.4 in all of the examples.

\subsection{One-Dimensional Examples}
\paragraph{Example 1---Shock-Entropy Wave Interaction Problem.} In the first example taken from \cite{Shu88}, we consider the shock-entropy
wave interaction problem. The initial conditions,
\begin{equation*}
(\rho, u,p)(x,0)=\begin{cases}(1.51695,0.523346,1.805),&x<-4.5,\\(1+0.1\sin(20x),0,1),&x>-4.5,\end{cases}
\end{equation*}
correspond to a forward-facing shock wave of Mach number 1.1 interacting with high-frequency density perturbations, that is, as the shock
wave moves, the perturbations spread ahead. We set the free boundary condition at the both ends of the computational domain $[-10,5]$.

\begin{figure}[ht!]
\centerline{\includegraphics[trim=0.9cm 0.4cm 1.2cm 0.6cm, clip, width=6.4cm]{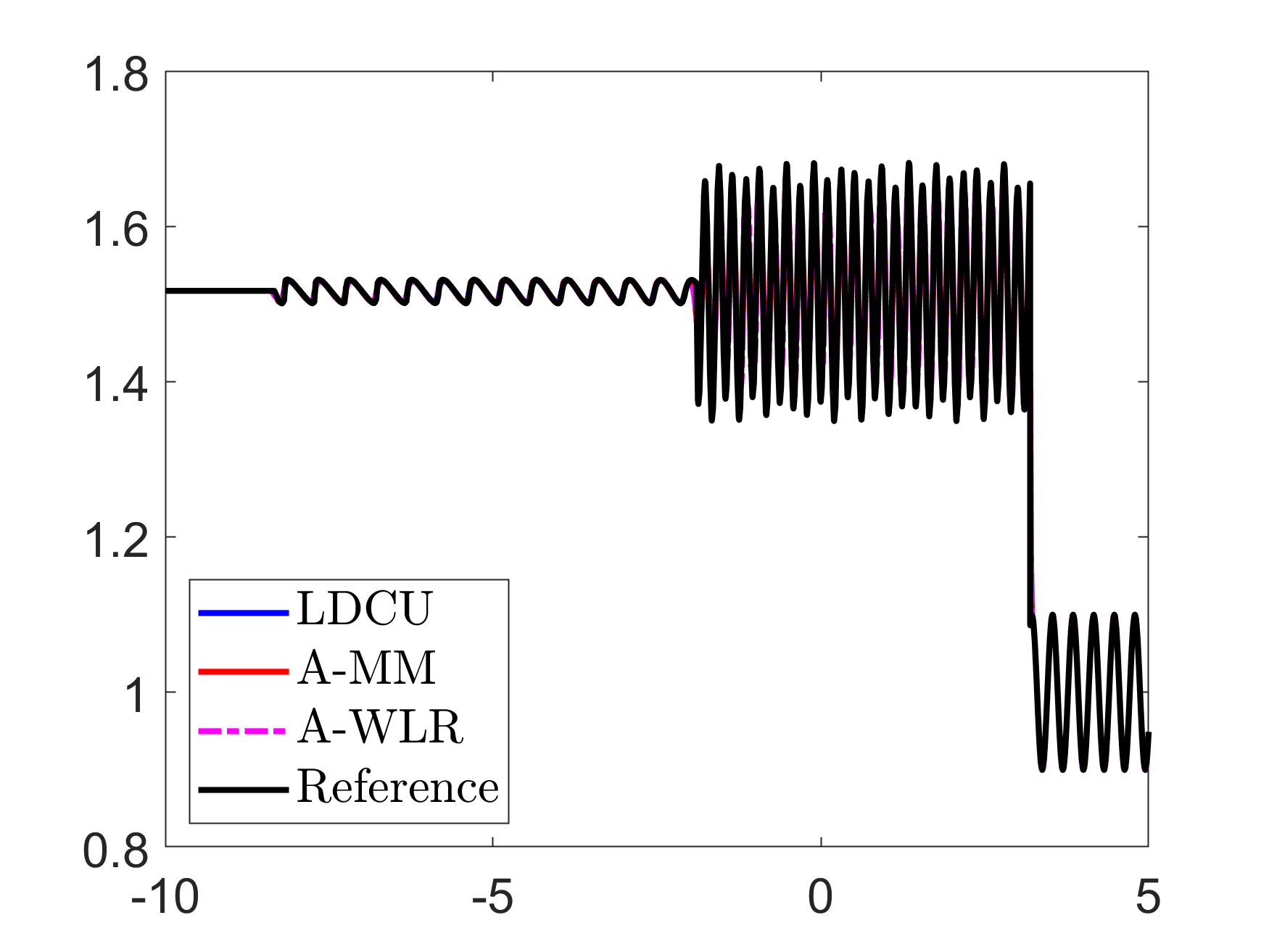}\hspace{1cm}
            \includegraphics[trim=0.9cm 0.4cm 1.2cm 0.6cm, clip, width=6.4cm]{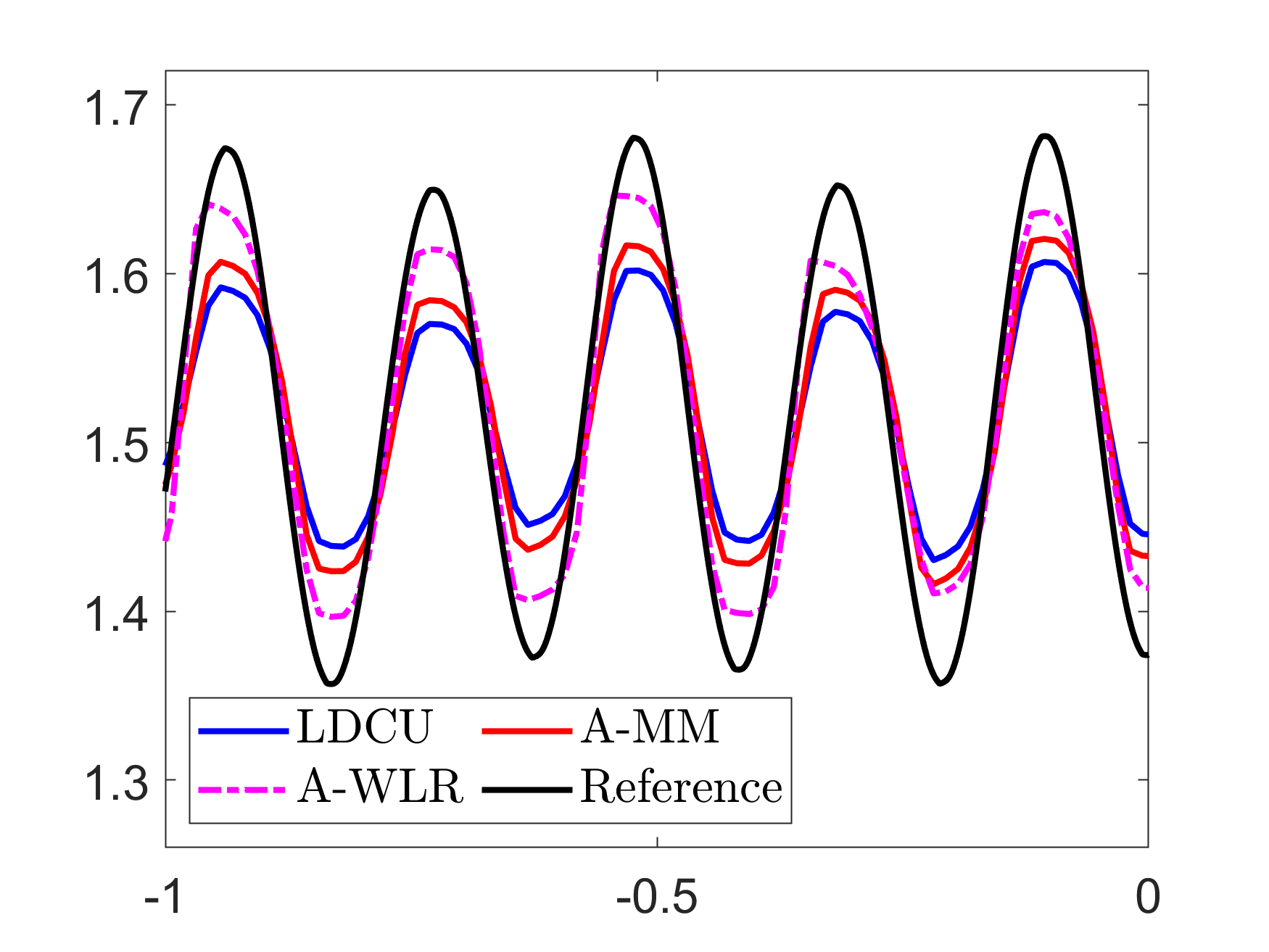}}
\caption{\sf Example 1: Density $\rho$ computed by the LDCU, A-MM, and A-WLR schemes (left) and zoom at $x\in[-1,0]$ (right).\label{fig3}}
\end{figure}

We compute the numerical solution until the final time $t=5$ by the LDCU, A-MM, and A-WLR schemes on a uniform mesh with $\dx=1/80$. We use
the adaption constant $\texttt{C}=0.1$ in the A-WLR scheme. The numerical results at time $t=5$ are presented in Figure \ref{fig3} along
with the reference solution computed by the LDCU scheme on a much finer mesh with $\dx=1/1600$. On the right panel of Figure \ref{fig3}, we
zoom the obtained solutions at the interval $[-1,0]$, at which the exact solution is smooth but has an oscillatory nature. As one can see,
this part of the solution is resolved much more accurately by the two adaptive schemes, especially by the A-WLR one.

\paragraph{Example 2---Shock-Density Wave Interaction Problem.} In the second example taken from \cite{SO89}, we consider the
shock-density wave interaction problem. The initial data,
\begin{equation*}
(\rho,u,p)(x,0)=\begin{cases}\Big(\dfrac{27}{7},\dfrac{4\sqrt{35}}{9},\dfrac{31}{3}\Big),&x<-4,\\[0.8ex](1+0.2\sin(5x),0,1),&x>-4,
\end{cases}
\end{equation*}
are prescribed in the computational domain $[-5,15]$ subject to the free boundary conditions.

\begin{figure}[ht!]
\centerline{\includegraphics[trim=0.9cm 0.4cm 1.1cm 0.6cm, clip, width=6.4cm]{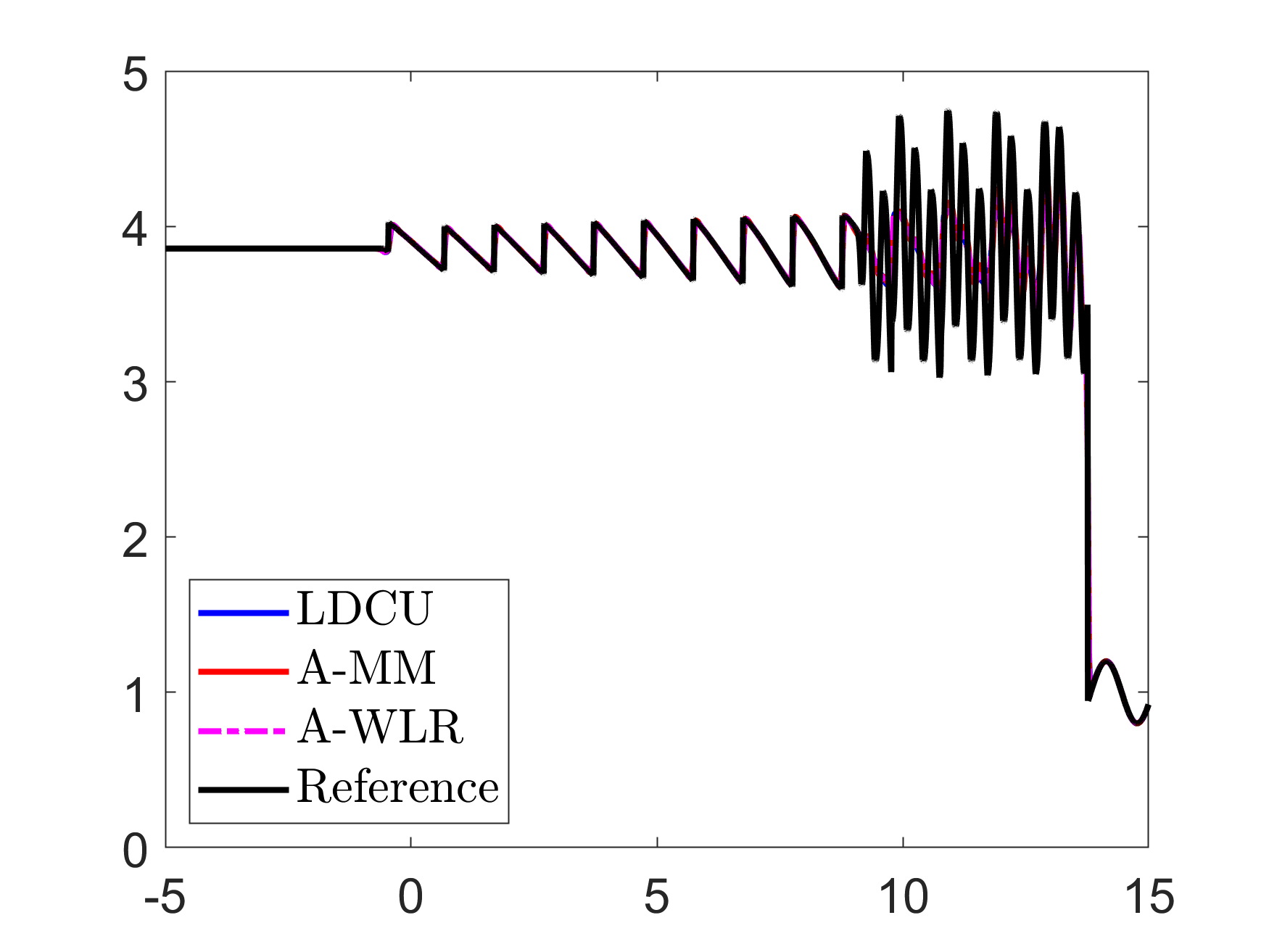}\hspace{1cm}
            \includegraphics[trim=0.9cm 0.4cm 1.1cm 0.6cm, clip, width=6.4cm]{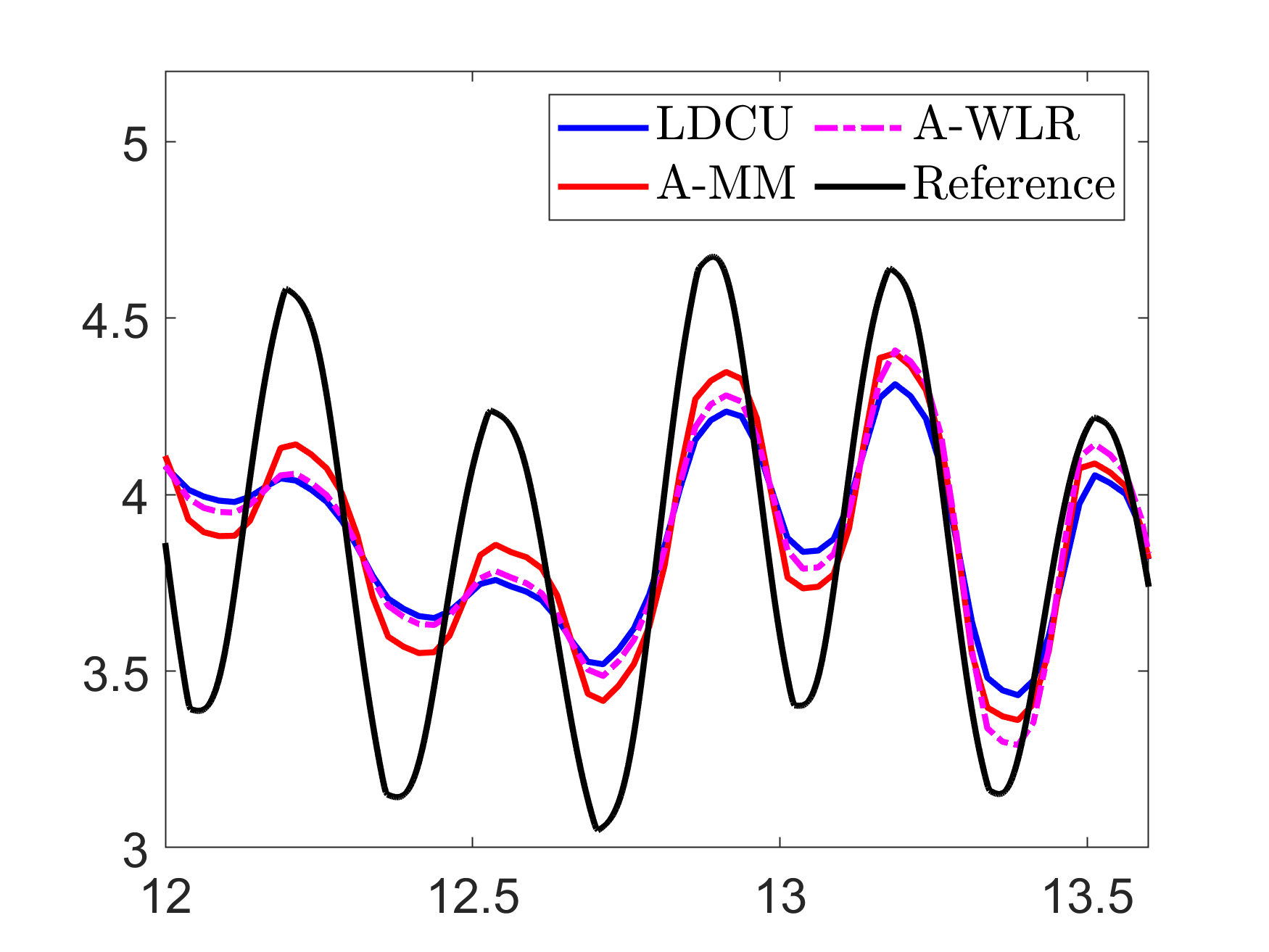}}
\caption{\sf Example 2: Density $\rho$ computed by the LDCU, A-MM, and A-WLR schemes (left) and zoom at $x\in[11.8,13.6]$ (right).
\label{fig4}}
\end{figure}

We compute the numerical solutions by the LDCU, A-MM, and A-WLR schemes on the uniform mesh with $\dx=1/40$ until the final time $t=5$. The
A-WLR scheme is used with the adaption constant $\texttt{C}=0.35$. We present the obtained numerical results in Figure \ref{fig4} together
with the reference solution computed by the LDCU scheme on a much finer mesh with $\dx=1/400$. It can be clearly seen in Figure \ref{fig4}
(right) that both of the adaptive schemes produce more accurate results compared to those obtained by the LDCU scheme. One can also observe
that in this example, unlike the previous one, the A-MM scheme achieves higher resolution of the smooth parts of the solution compared with
its A-WLR counterpart. This is attributed to a relatively large value of $\texttt{C}$ used in this example.

\begin{figure}[ht!]
\centerline{\includegraphics[trim=0.9cm 0.4cm 1.1cm 0.2cm, clip, width=6.4cm]{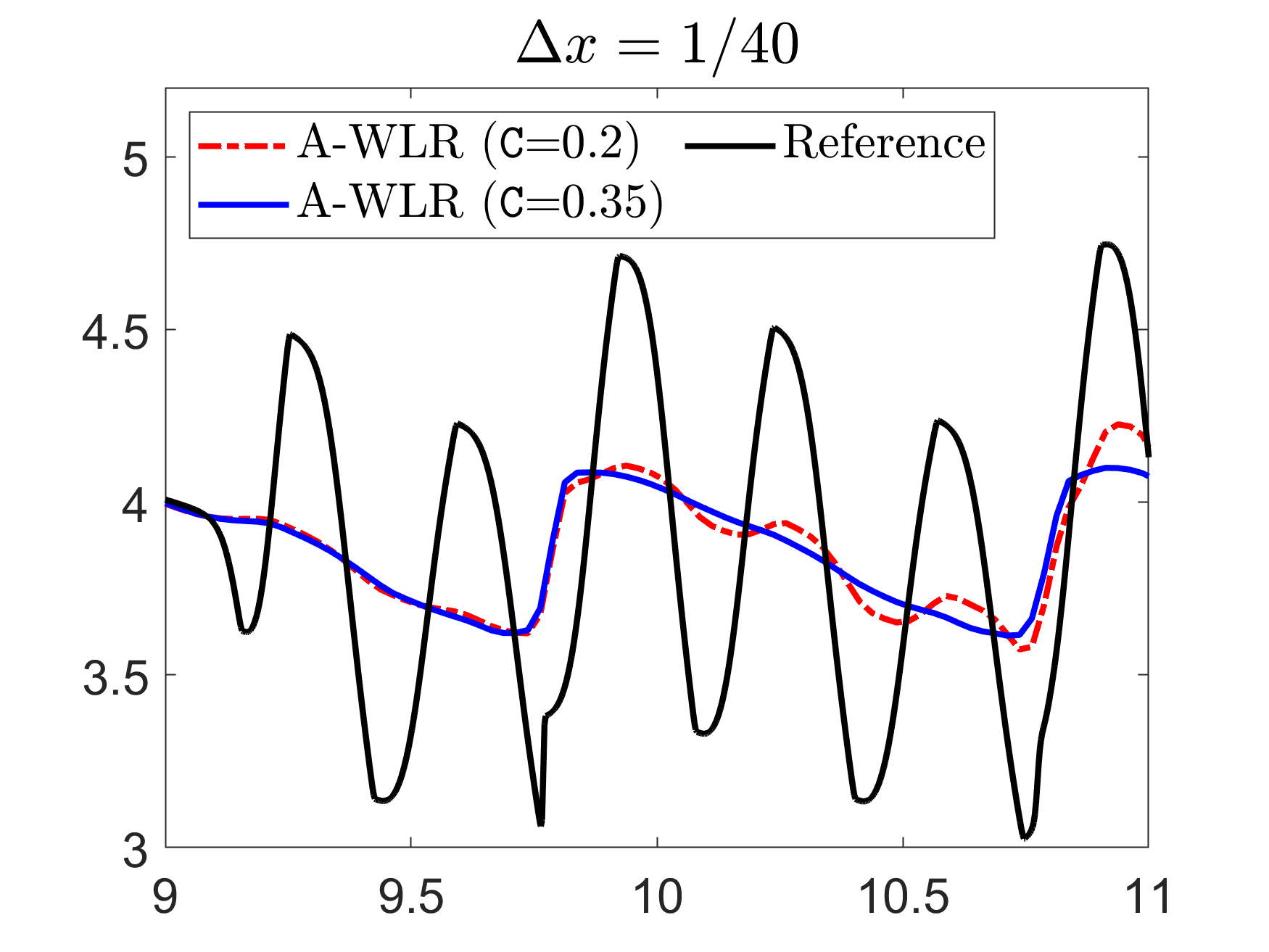}\hspace{1cm}
            \includegraphics[trim=0.9cm 0.4cm 1.1cm 0.2cm, clip, width=6.4cm]{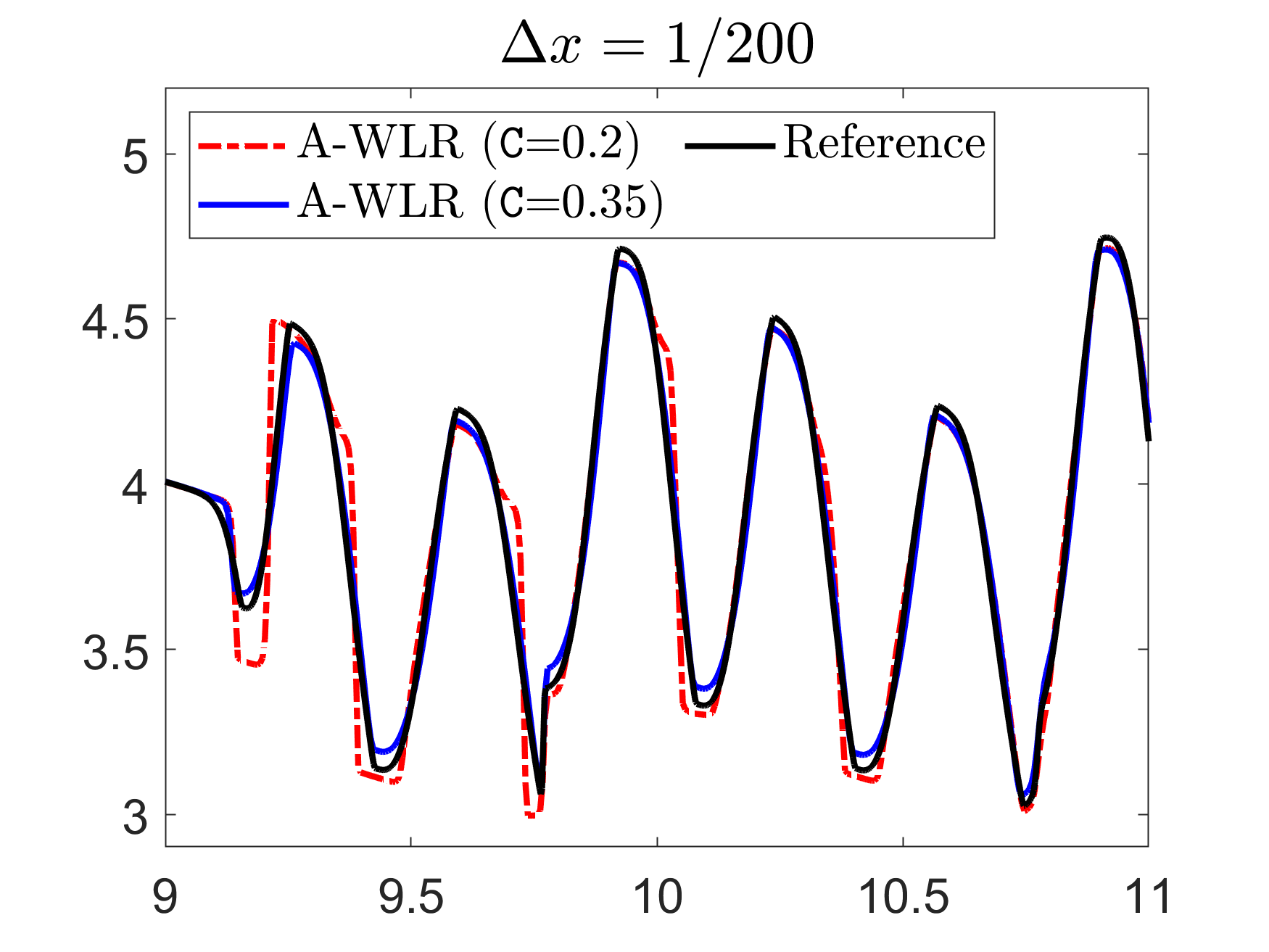}}
\caption{\sf Example 2: Density $\rho$ computed by the A-WLR scheme with $\texttt{C}=0.2$ and $\texttt{0.35}$ on a coarse (left) and fine
(right) meshes. Zoom at $x\in[9,11]$.\label{fig4a}}
\end{figure}

Recall that one of the key points in the proposed WLR-based scheme adaption strategy is tuning the adaption constant $\texttt{C}$. In
\cite{Kurganov12a}, where a WLR-based adaptive artificial viscosity was introduced and studied, the viscosity coefficient, which is directly
related to $\texttt{C}$, was first adjusted on a coarse mesh and then used for the high-resolution computations on finer meshes. However,
this strategy does not seem to be robust in the A-WLR scheme as the numerical results computed by the A-WLR scheme may still have
staircase-like structure in the areas where the coarse mesh solution is smooth. In order to illustrate this, we use a slightly smaller
adaption constant $\texttt{C}=0.2$ and compute the results with $\texttt{C}=0.2$ and $\texttt{C}=0.35$ on the uniform meshes with
$\dx=1/40$ and $1/200$. We present the obtained densities (zoomed at $x\in[9,11]$) in Figure \ref{fig4a}, where one can see that even though
the solution computed on the coarse mesh with $\texttt{C}=0.2$ is smooth, it develops clear staircase-like structures when the mesh is
refined.

\paragraph{Example 3---Blast Wave Problem.}
In the last 1-D example, we consider the strong shocks interaction problem from \cite{Woodward88}, which is considered on the interval
$[0,1]$ with the solid wall boundary conditions at both ends and subject to the following initial conditions:
\begin{equation*}
(\rho,u,p)(x,0)=\begin{cases}(1,0,1000),&x<0.1,\\(1,0,0.01),&0.1\le x\le 0.9,\\(1,0,100),&x>0.9.\end{cases}
\end{equation*}

\begin{figure}[ht!]
\centerline{\includegraphics[trim=1.3cm 0.4cm 0.8cm 0.6cm, clip, width=6.4cm]{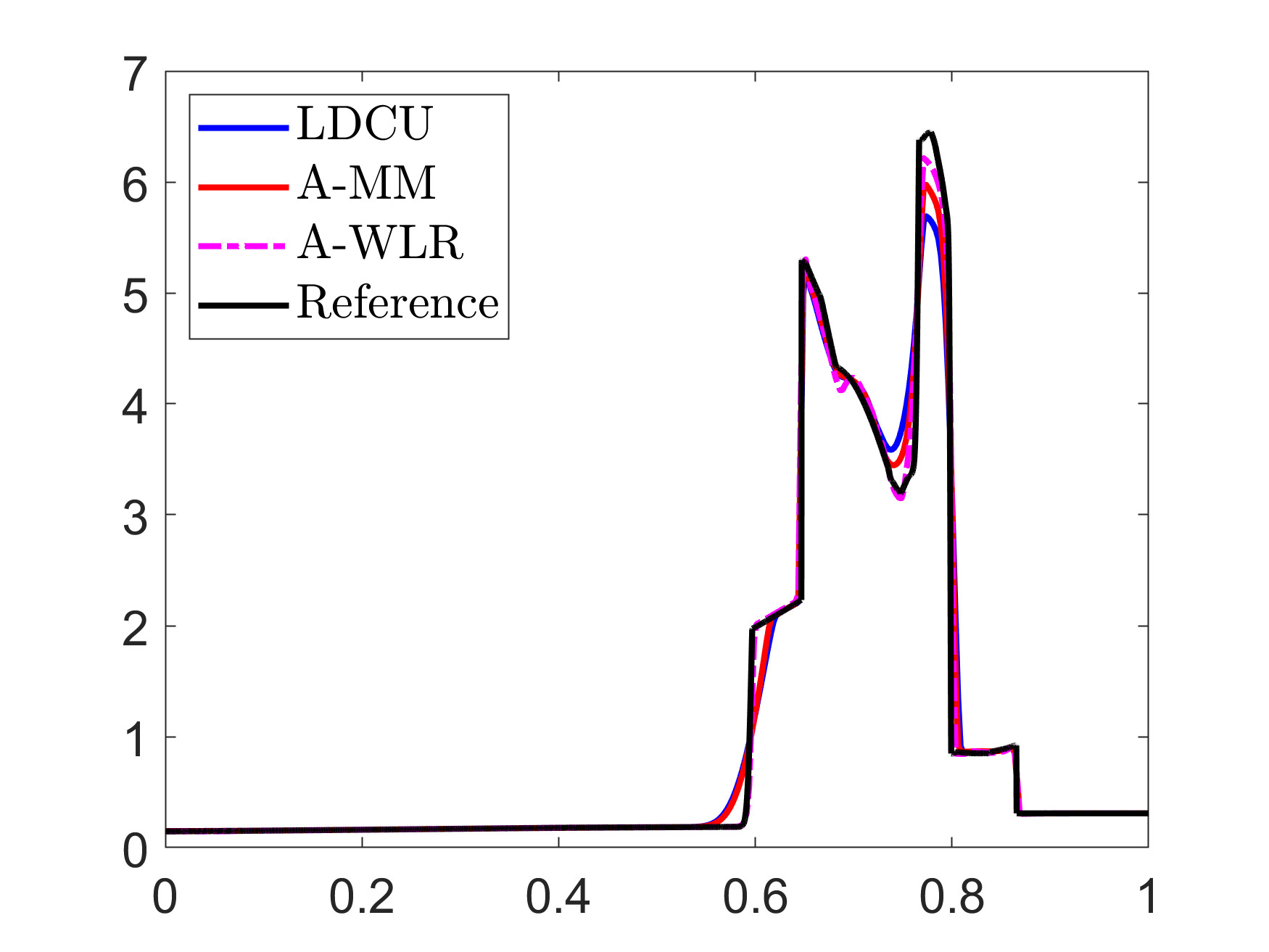}\hspace*{1.0cm}
            \includegraphics[trim=1.3cm 0.4cm 0.8cm 0.6cm, clip, width=6.4cm]{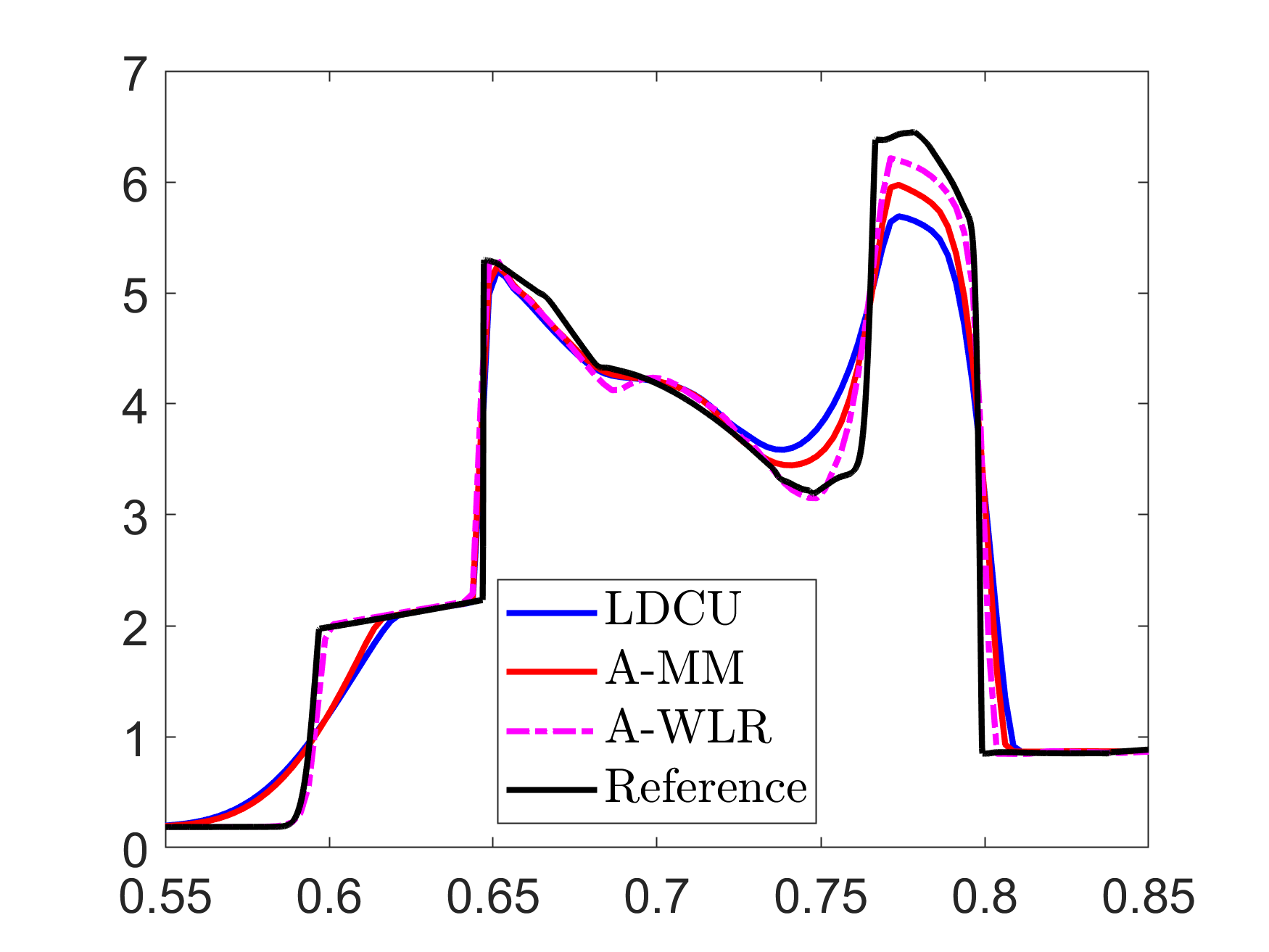}}
\caption{\sf Example 3: Density $\rho$ computed by the LDCU, A-MM, and A-WLR schemes (left) and zoom at $x\in[0.55,0.85]$.\label{fig5}}
\end{figure}

We compute the numerical solutions until the final time $t=0.038$ by the LDCU, A-MM, and A-WLR (with the adaption constant $\texttt{C}=0.1$)
schemes on a uniform mesh with $\dx=1/400$ and implement the LDCU scheme on a much finer grid with $\dx=1/4000$ to compute the reference
solution. The obtained results, presented in Figure \ref{fig5}, demonstrate that while the A-MM scheme outperforms the LDCU one, the A-WLR
results are even more accurate and the A-WLR scheme is capable of achieving a superb resolution of the contact wave located at about
$x=0.6$. It is well-known that this contact wave is the one, which is hardest to get sharply resolved, and the A-WLR scheme is, to best of
our knowledge, the first Riemann-problem solver-free scheme that can achieve this goal.

\subsection{Two-Dimensional Examples}
\paragraph{Example 4---2-D Riemann Problem.} In the first 2-D example, we consider Configuration 3 of the 2-D Riemann problems from
\cite{Kurganov02} (see also \cite{Schulz93,Schulz93a,Zheng01}) with the initial conditions
\begin{equation*}
(\rho,u,v,p)(x,y,0)=\begin{cases}(1.5,0,0,1.5),&x>1,~y>1,\\(0.5323,1.206,0,0.3),&x<1,~y>1,\\(0.138,1.206,1.206,0.029),&x<1,~y<1,\\
(0.5323,0,1.206,0.3),&x>1,~y<1,\end{cases}
\end{equation*}
prescribed in the computational domain $[0,1.2]\times[0,1.2]$ subject to the free boundary conditions.

\begin{figure}[ht!]
\centerline{\includegraphics[trim=2.3cm 0.5cm 2.3cm 0.1cm, clip, width=5.6cm]{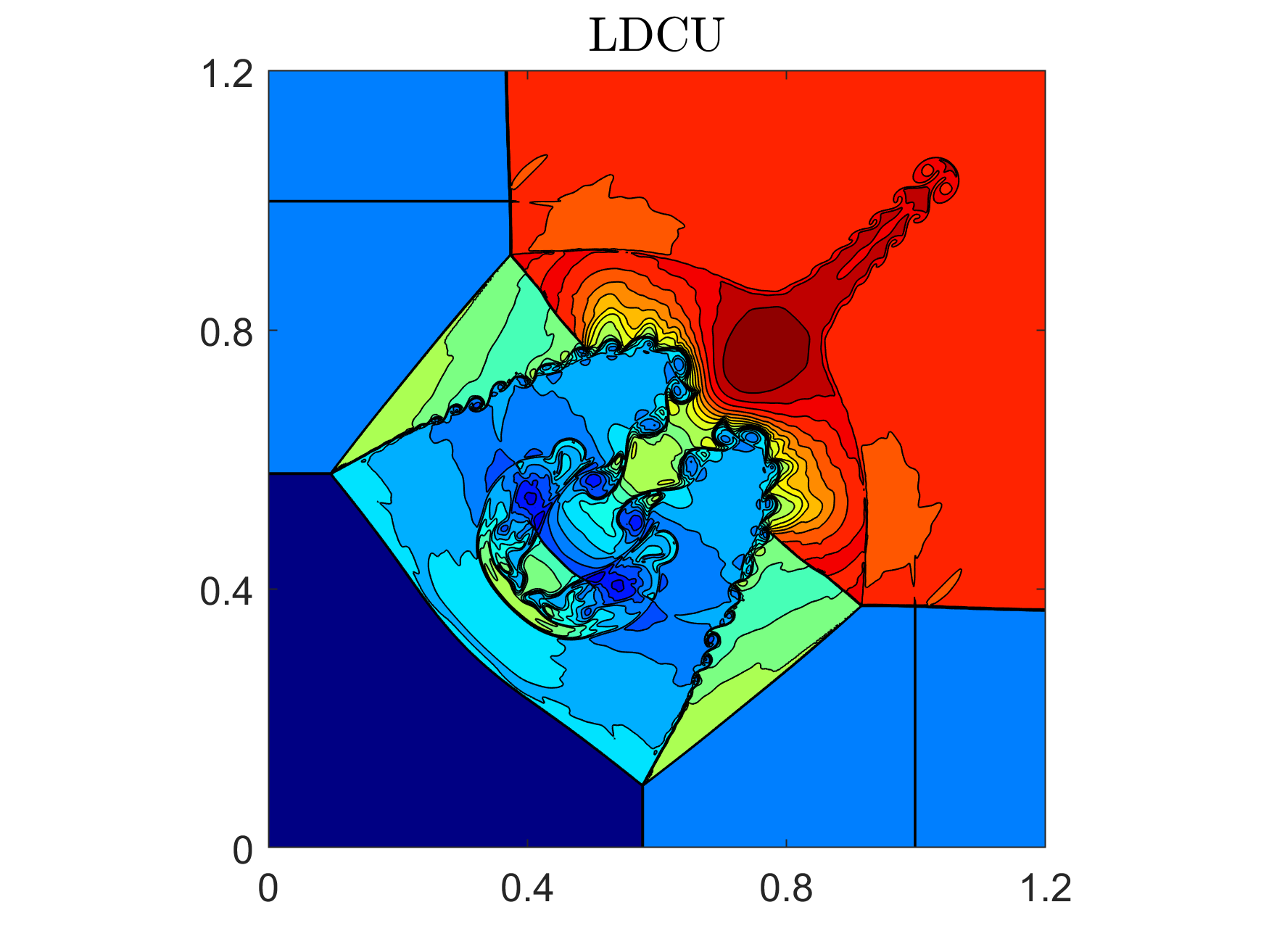}\hspace*{-0.0cm}
            \includegraphics[trim=3.0cm 0.5cm 2.3cm 0.1cm, clip, width=5.2cm]{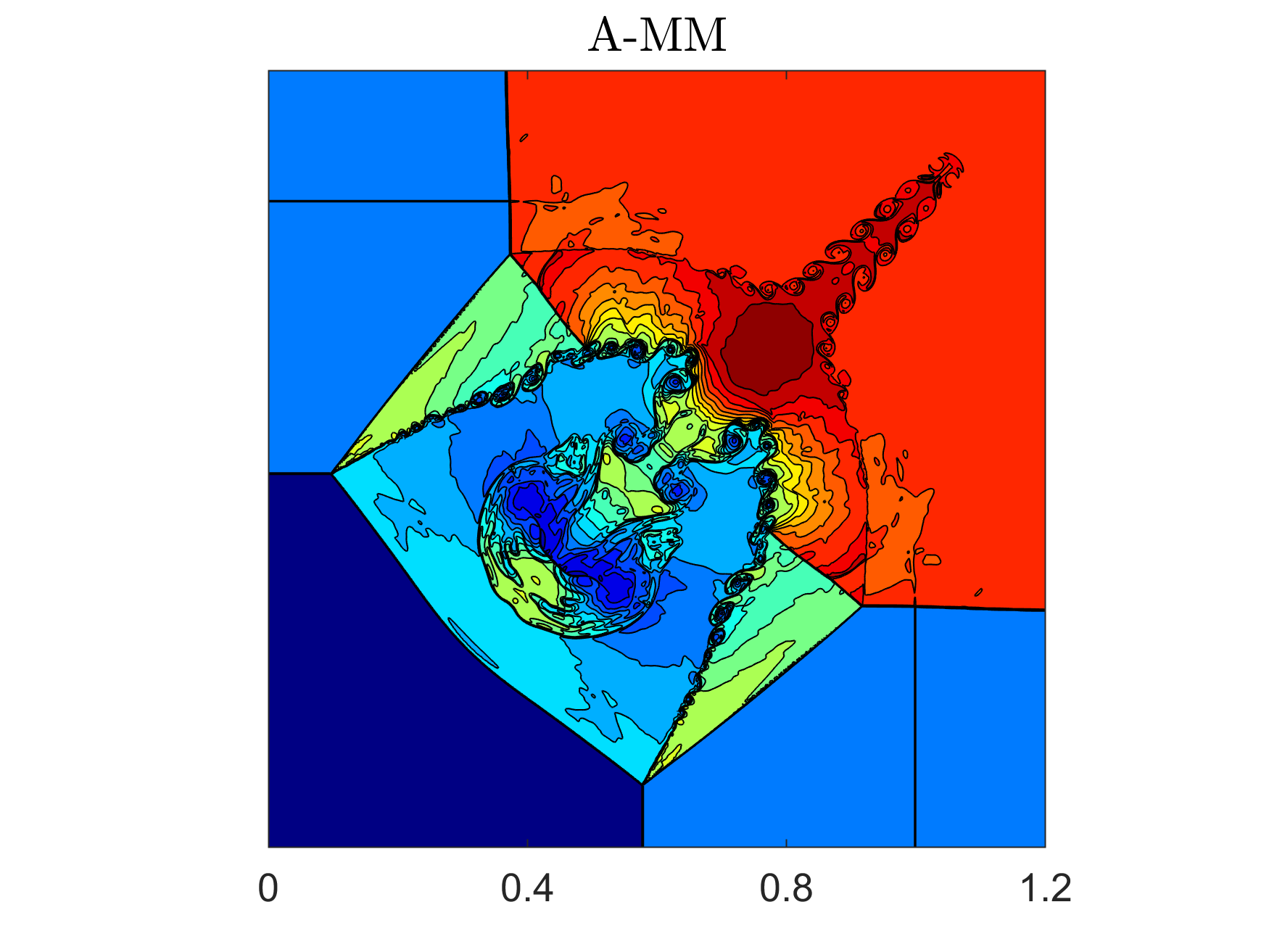}\hspace*{-0.0cm}
            \includegraphics[trim=1.9cm 0.5cm 1.6cm 0.1cm, clip, width=6.2cm]{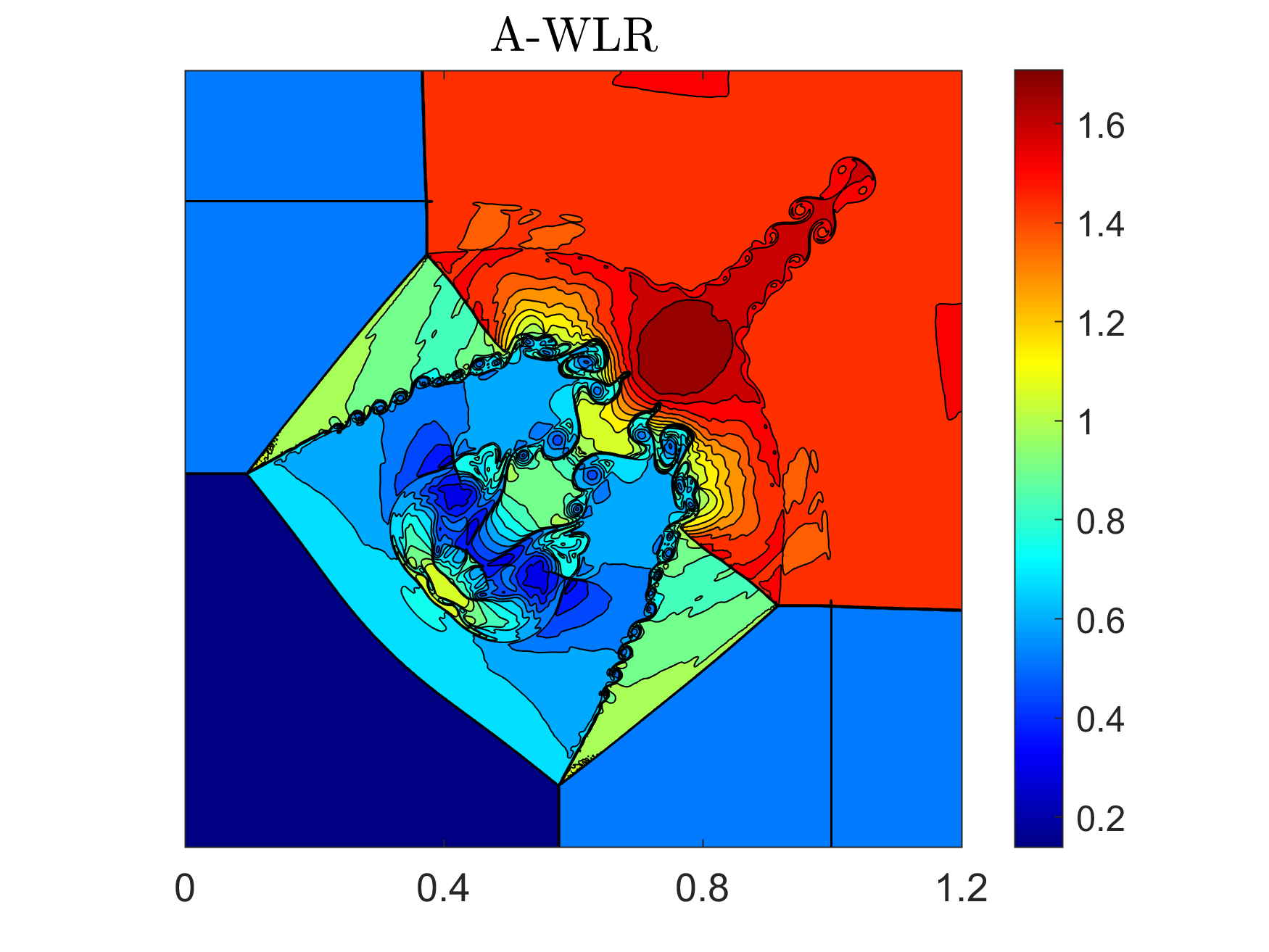}}
\caption{\sf Example 4: Density $\rho$ computed by the LDCU (left), A-MM (middle), and A-WLR (right) schemes.\label{fig6a}}
\end{figure}

We compute the numerical solution until the final time $t=1$ by the LDCU, A-MM, and A-WLR (with the adaption constant $\texttt{C}=4$)
schemes on the uniform mesh with $\dx=\dy=3/2500$ and present the obtained results in Figure \ref{fig6a}, where one can see that both the
A-MM and A-WLR schemes outperform the LDCU scheme in capturing a sideband instability of the jet in the zones of strong along-jet velocity
shear and the instability along the jets neck.

\begin{figure}[ht!]
\centerline{\includegraphics[trim=2.3cm 0.5cm 2.3cm 0.1cm, clip, width=5.6cm]{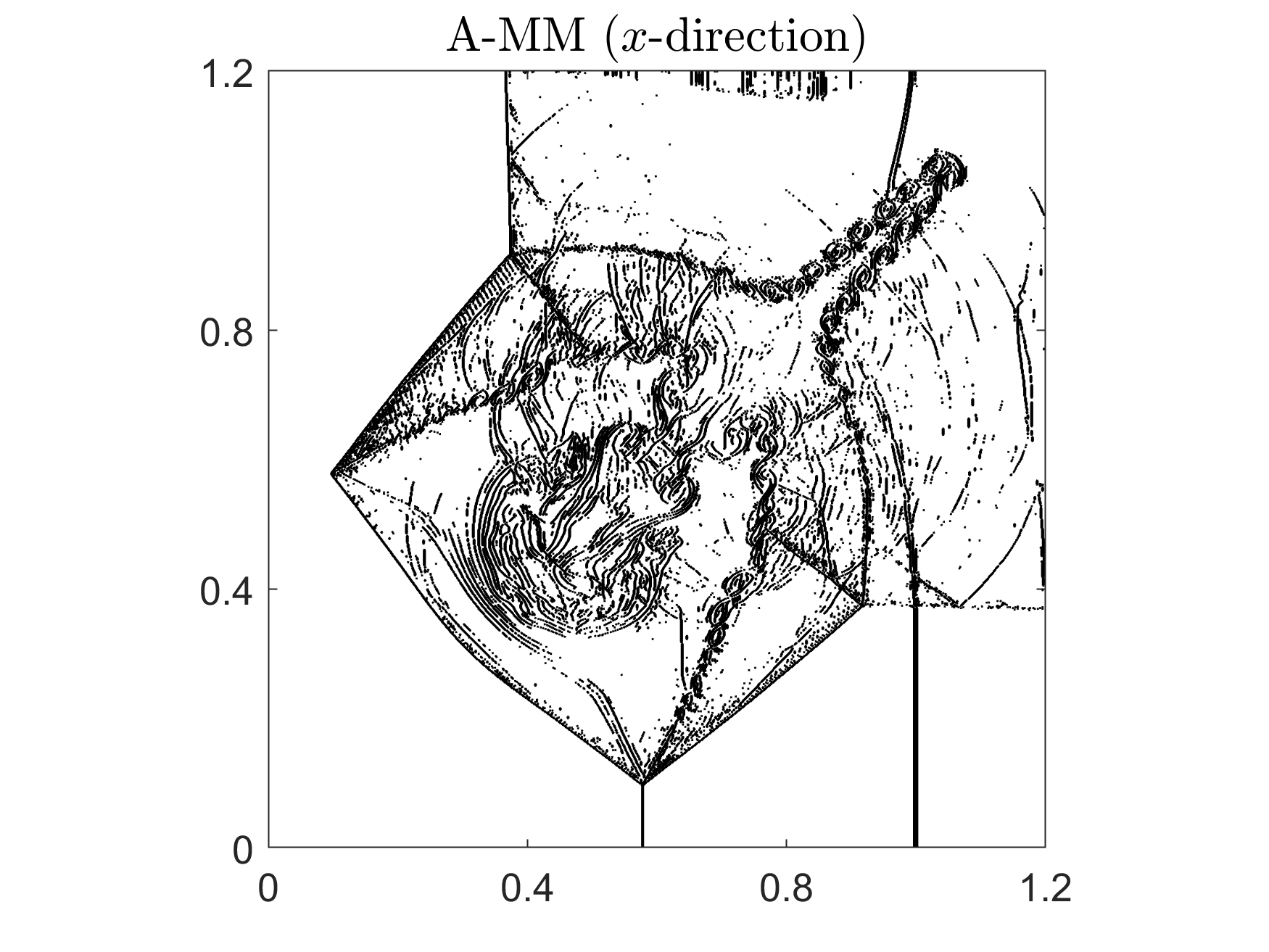}\hspace*{0.2cm}
            \includegraphics[trim=2.9cm 0.5cm 2.3cm 0.1cm, clip, width=5.25cm]{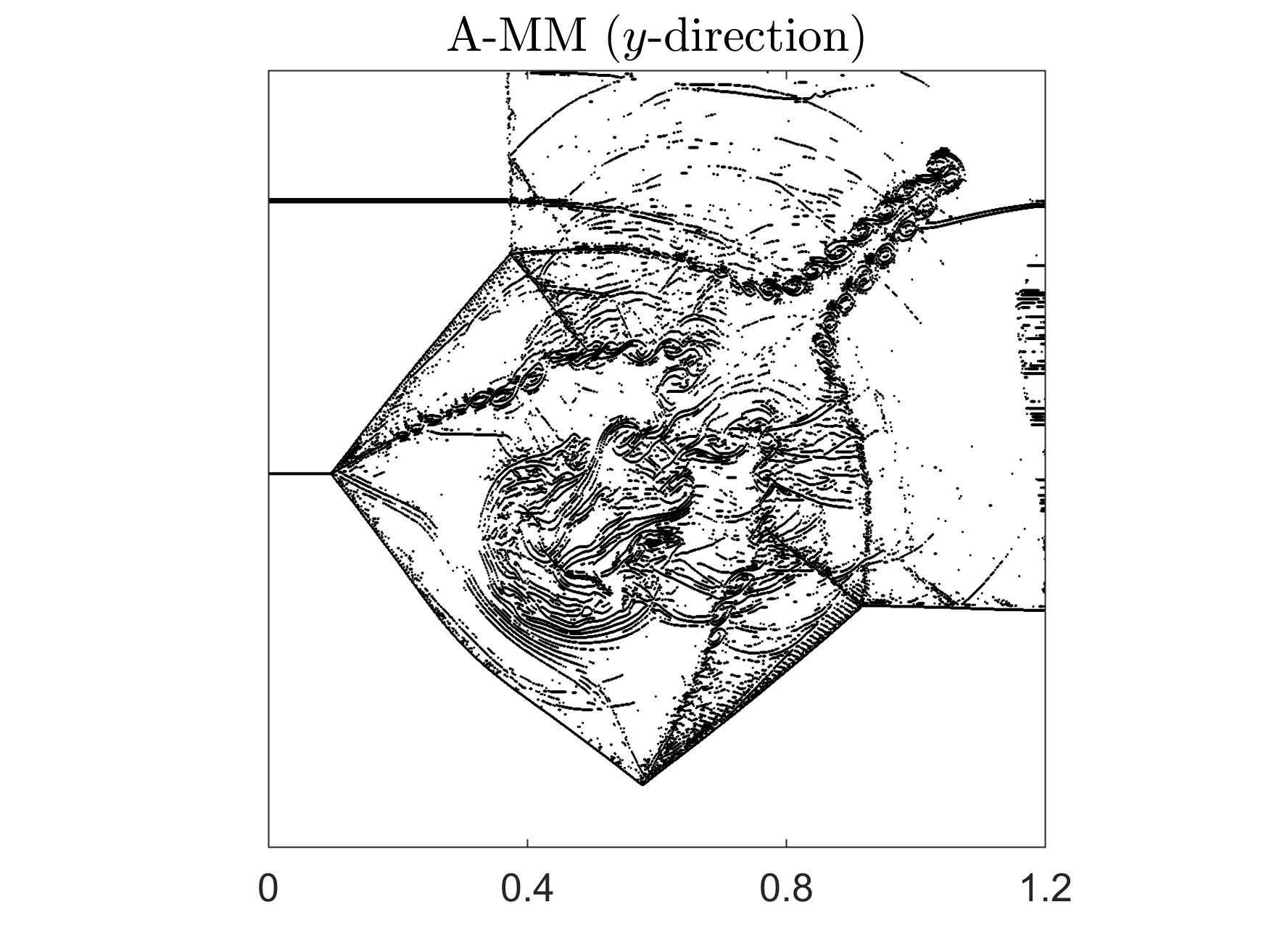}\hspace*{0.5cm}
            \includegraphics[trim=2.9cm 0.5cm 2.3cm 0.1cm, clip, width=5.25cm]{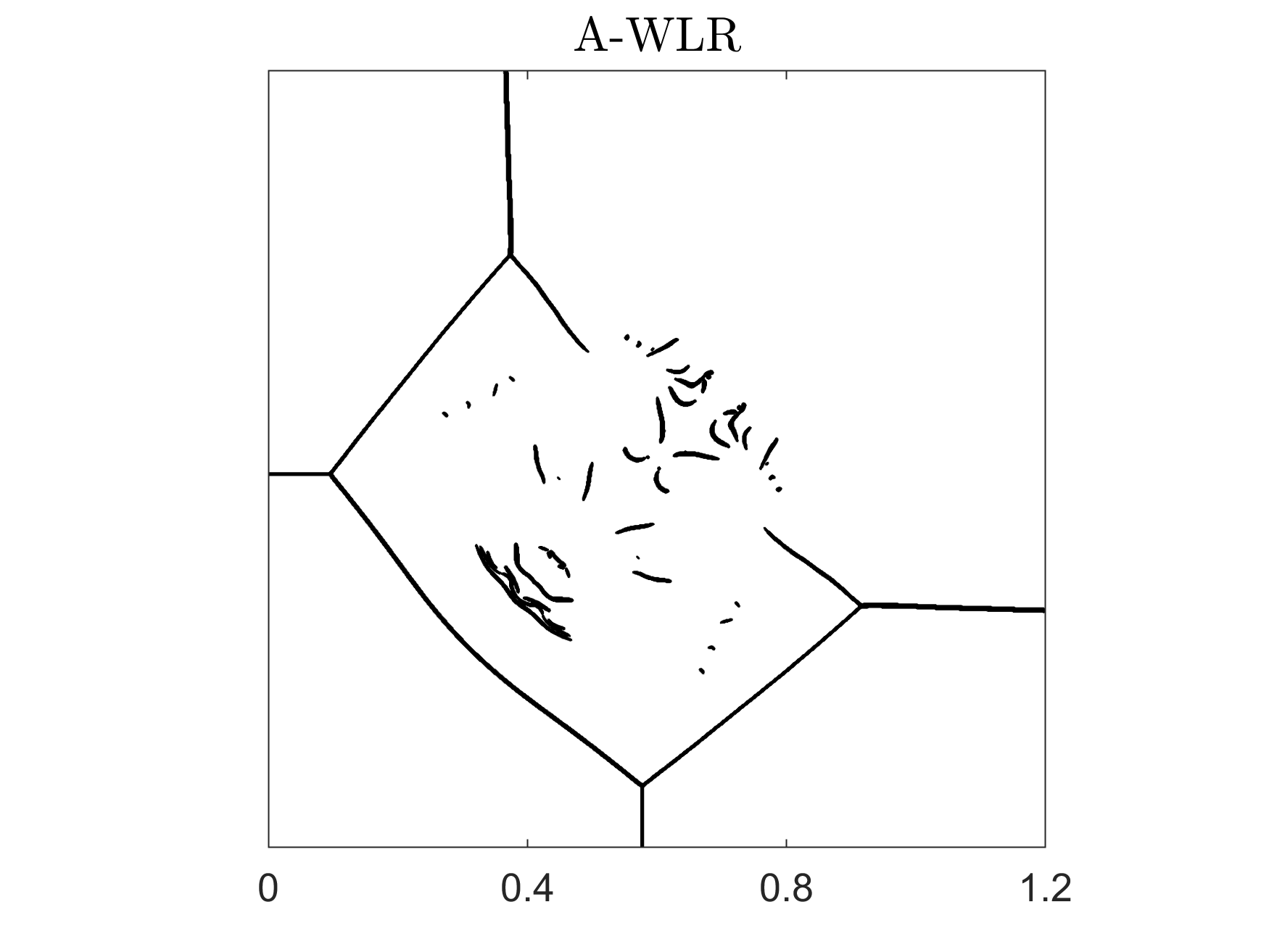}}
\caption{\sf Example 4: Areas detected as having large $x$- (left) and $y$-derivatives (middle) by the MM-based SI and the ``rough'' areas
detected by the WLR-based SI (right).\label{fig6b}}
\end{figure}

In Figure \ref{fig6b}, we show the regions which have been detected as ``rough'' by the SIs at the final time. We first indicate (in the
left and middle panels) the regions in which the MM-based SI detected large $x$- and $y$-directional derivatives, respectively, and where
the overcompressive SBM directional limiters have been used. In the right panel, we show the ``rough'' regions indicated by the WLR-based
SI. As one can see, when the A-WLR scheme is used, a sharper SBM limiter is implemented only in a small part of the computational domain,
mostly around the shocks.

\paragraph{Example 5---Implosion Problem.} In this example, we consider the implosion problem taken from \cite{Liska03}. The initial
conditions,
\begin{equation*}
(\rho,u,v,p)(x,y,0)=\begin{cases}(0.125,0,0,0.14),&|x|+|y|<0.15,\\(1,0,0,1),&\mbox{otherwise},\end{cases}`
\end{equation*}
are prescribed in $[0,0.5]\times[0,0.5]$ with the solid wall boundary conditions imposed at all of the four sides. This example was designed
to test the amount of numerical diffusion present in different schemes as there is a jet forming near the origin and propagating along the
diagonal $y=x$ direction, and schemes containing large numerical diffusion may not resolve the jet at all or the jet propagation velocity
may be affected by the numerical diffusion.

\begin{figure}[ht!]
\centerline{\includegraphics[trim=2.3cm 0.5cm 2.3cm 0.1cm, clip, width=5.6cm]{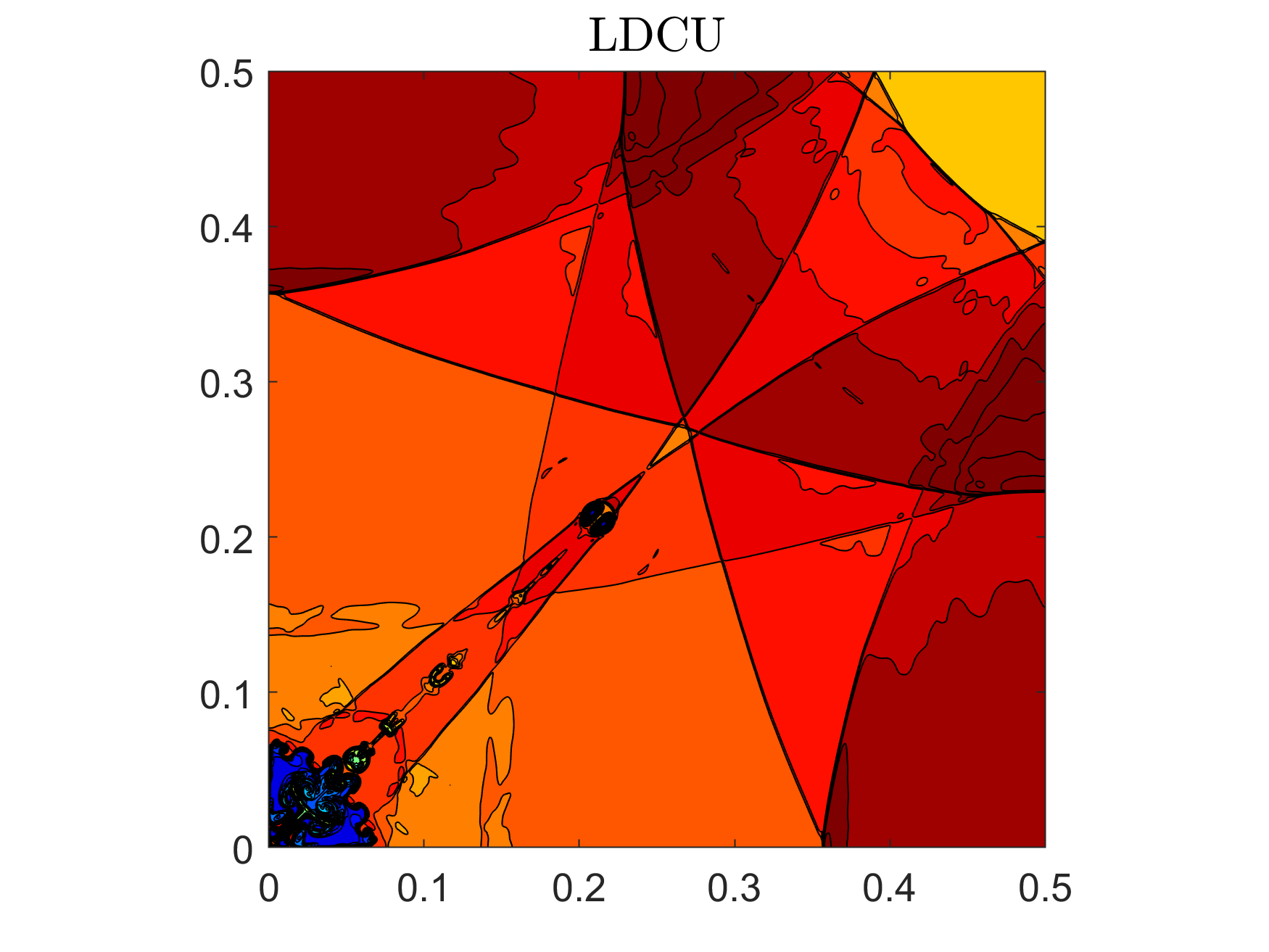}\hspace*{-0.0cm}
            \includegraphics[trim=3.0cm 0.5cm 2.3cm 0.1cm, clip, width=5.2cm]{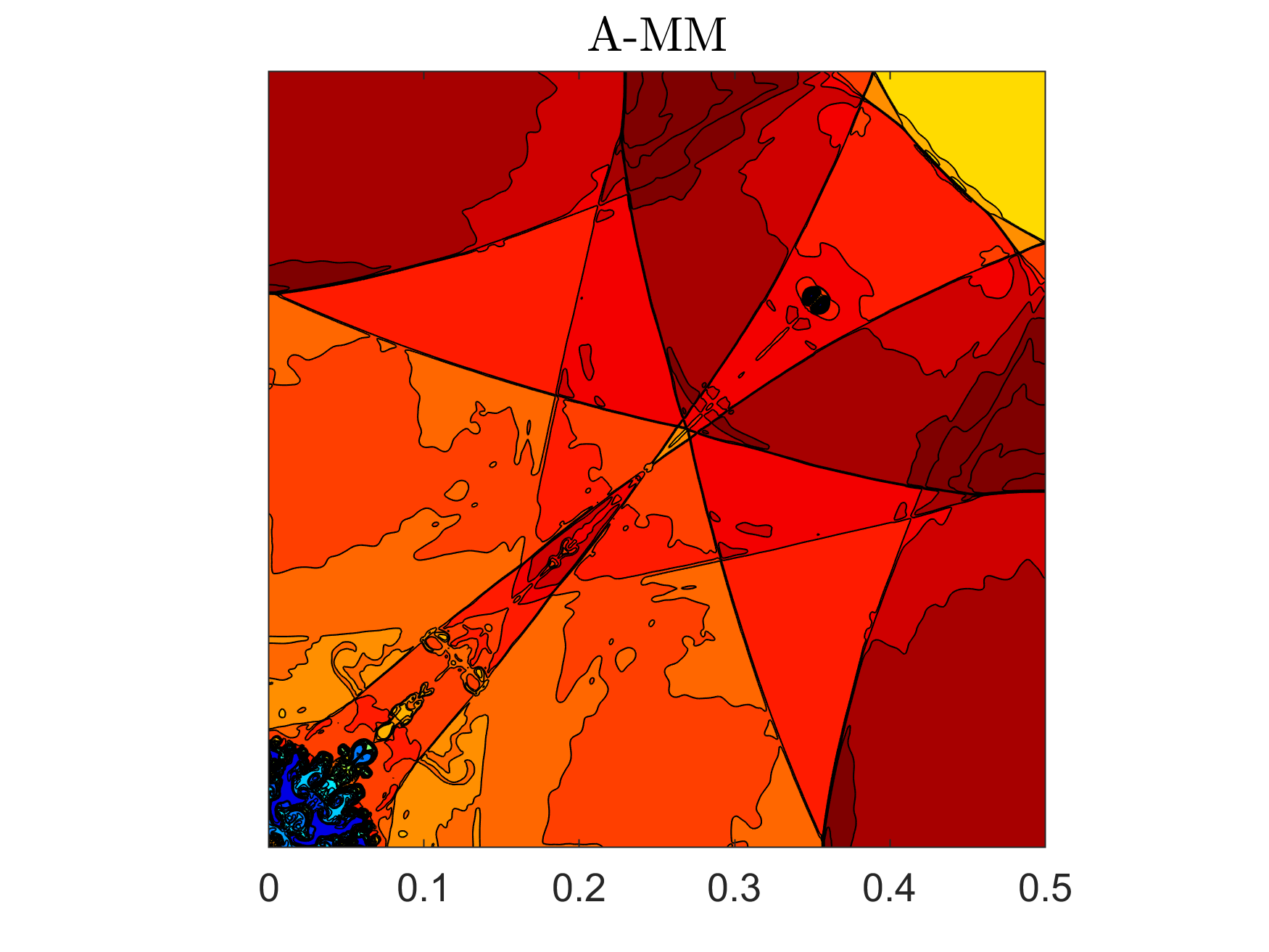}\hspace*{-0.0cm}
            \includegraphics[trim=1.9cm 0.5cm 1.6cm 0.1cm, clip, width=6.2cm]{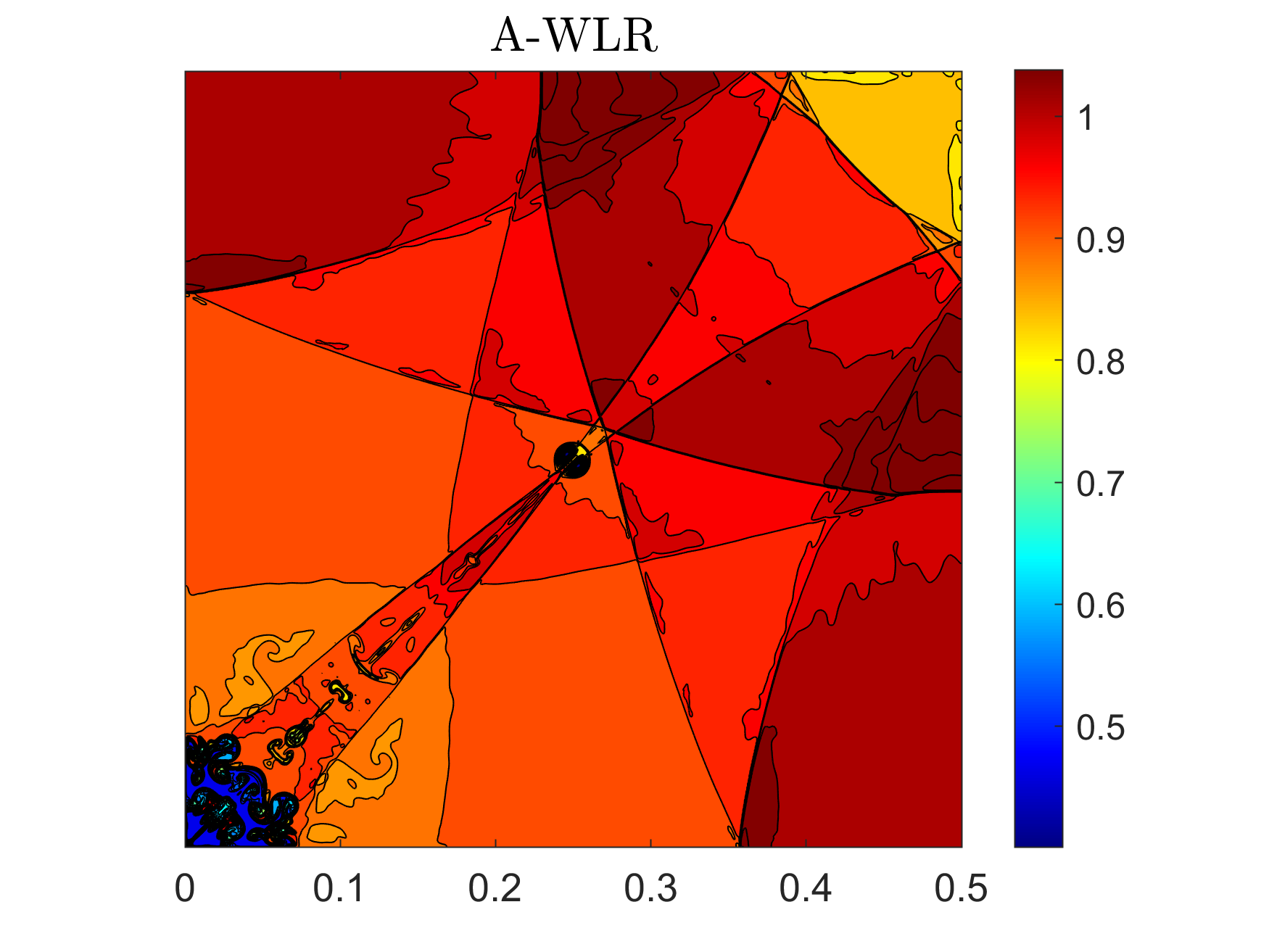}}
\caption{\sf Example 5: Density $\rho$ computed by the LDCU (left), A-MM (middle), and A-WLR (right) schemes.\label{fig8a}}
\end{figure}

\begin{figure}[ht!]
\centerline{\includegraphics[trim=2.3cm 0.5cm 2.3cm 0.1cm, clip, width=5.6cm]{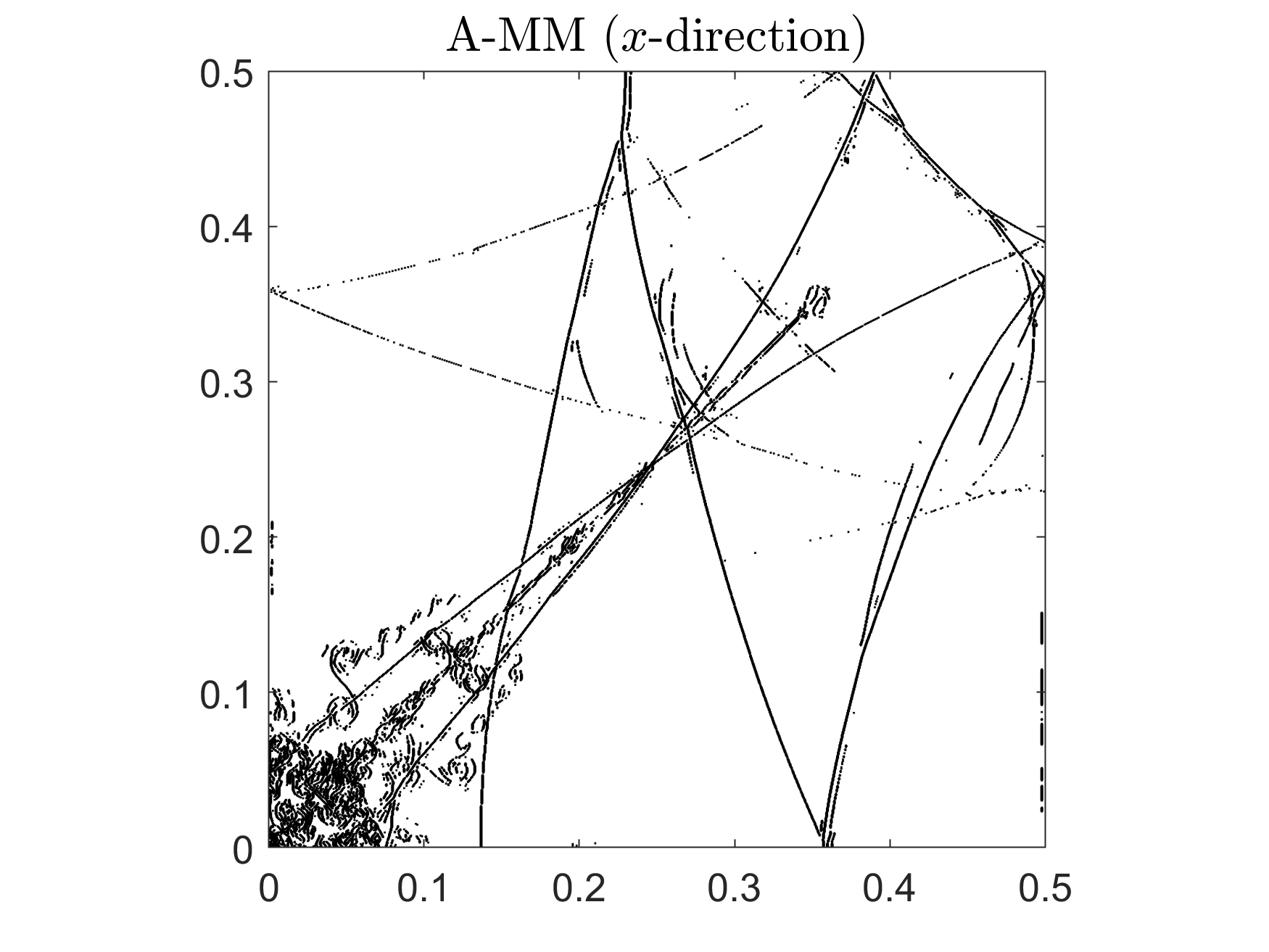}\hspace*{0.2cm}
            \includegraphics[trim=2.9cm 0.5cm 2.3cm 0.1cm, clip, width=5.25cm]{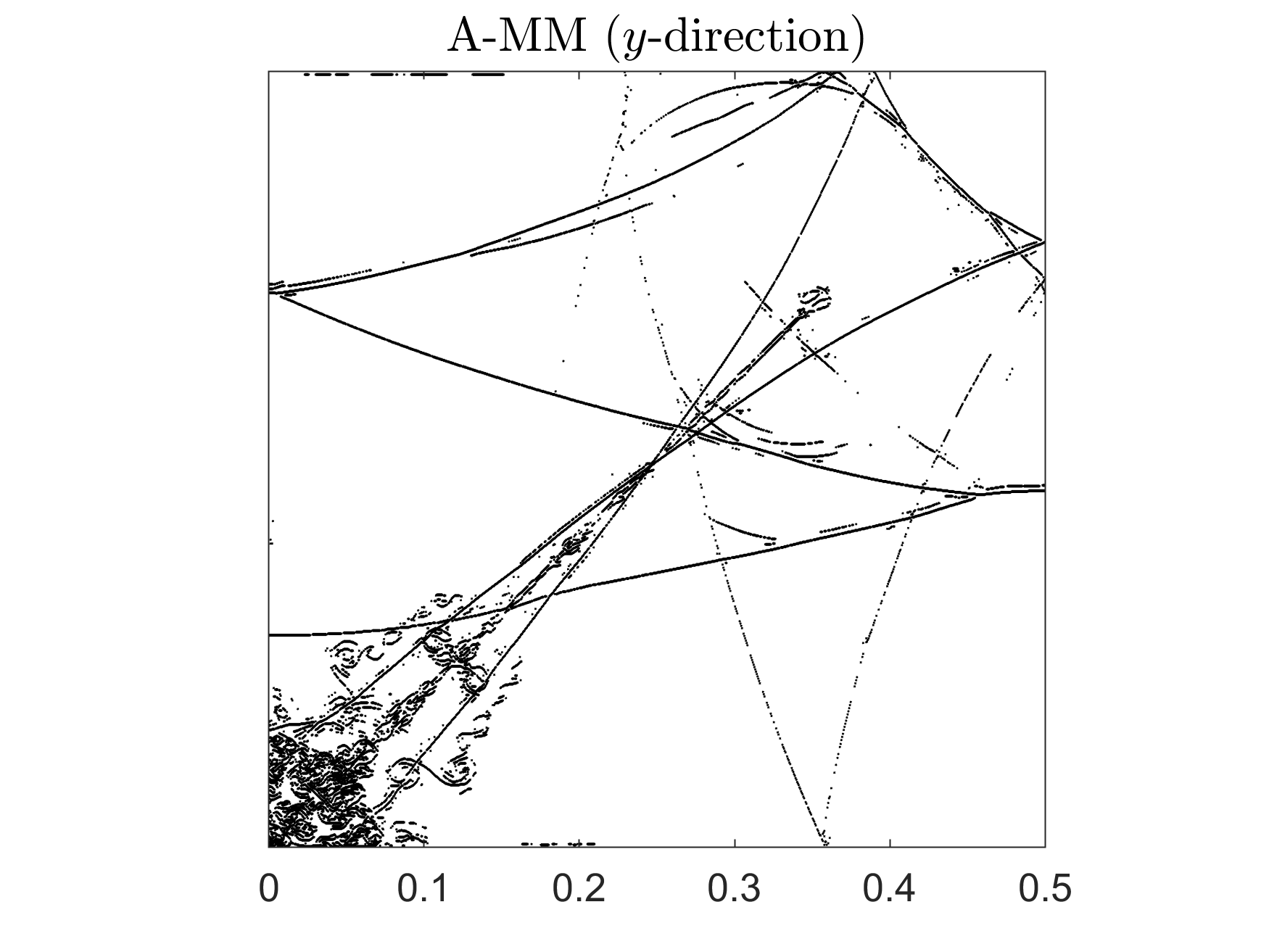}\hspace*{0.5cm}
            \includegraphics[trim=2.9cm 0.5cm 2.3cm 0.1cm, clip, width=5.25cm]{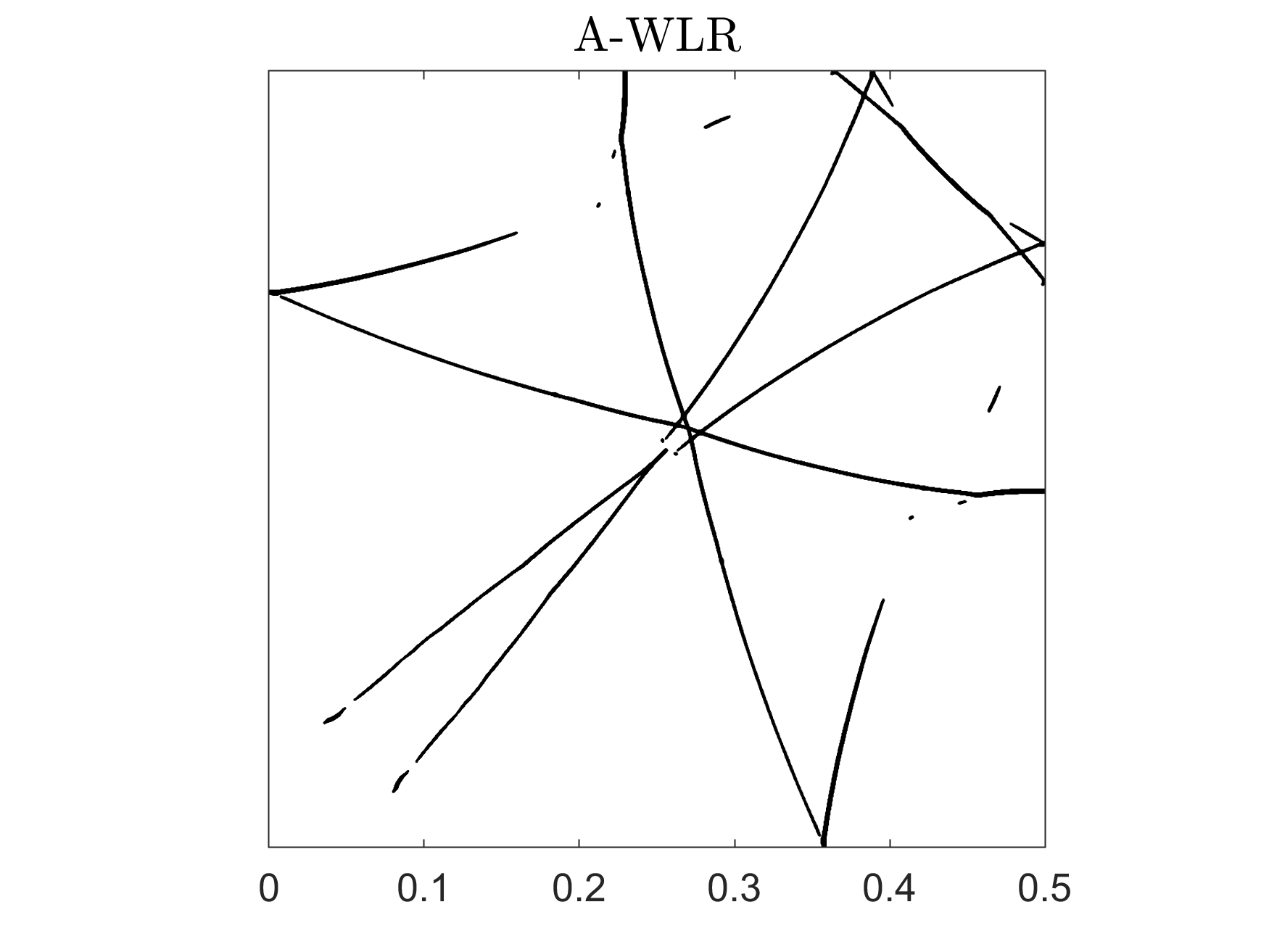}}
\caption{\sf Example 5: Areas detected as having large $x$- (left) and $y$-derivatives (middle) by the MM-based SI and the ``rough'' areas
detected by the WLR-based SI (right).\label{fig8b}}
\end{figure}

We compute the numerical solution until the final time $t=2.5$ by the LDCU, A-MM, and A-WLR (with the adaption constant $\texttt{C}=5$)
schemes on the uniform mesh with $\dx=\dy=1/2000$ and present the obtained numerical results in Figure \ref{fig8a}. As one can observe,
while the jet is generated by all of the studied schemes, it propagates much further in the diagonal direction when the solution is computed
by one of the adaptive schemes. In fact, the A-MM scheme seems to contain even smaller amount of numerical diffusion than the A-WLR one.
This can be confirmed by the results presented in Figure \ref{fig8b}, where we see that the ``rough'' areas detected at the final time by
the A-WLR scheme are concentrated along the shock waves only. This means that the A-WLR scheme uses the overcompressive SBM limiter in a
smaller part of the computational domain compared with the A-MM scheme. It is worth noting that the A-WLR scheme can be made less
dissipative by decreasing the adaption constant $\texttt{C}$, but it might be difficult to tune $\texttt{C}$ in this example as for smaller
values of $\texttt{C}$ the A-WLR scheme may produce reasonably sharp results on a coarse mesh, but it may develop severe instabilities when
the mesh is refined to $\dx=\dy=1/2000$.
\begin{rmk}
The solution of the studied initial-boundary value problem is symmetric with respect to the axis $y=x$, but this symmetry may be destroyed
by the roundoff errors when the solution is computed by the studied low-dissipation schemes. In order to prevent the loss of symmetry, we
have used a very simple strategy introduced in \cite{WDGK2020}: upon completion of each time evolution step, we replace the computed values
$\xbar\mU_{j,\,k}$ with $\widehat{\bm U}_{j,k}$, where
\begin{equation*}
\widehat\rho_{j,k}:=\frac{\xbar\rho_{j,k}+\xbar\rho_{k,j}}{2},~~
(\widehat{\rho u})_{j,k}:=\frac{(\xbar{\rho u})_{j,k}+(\xbar{\rho v})_{k,j}}{2},~~
(\widehat{\rho v})_{j,k}:=\frac{(\xbar{\rho v})_{j,k}+(\xbar{\rho u})_{k,j}}{2},~~
\widehat E_{j,k}:=\frac{\xbar E_{j,k}+\xbar E_{k,j}}{2},
\end{equation*}
for all $j,k$. For more sophisticated symmetry enforcement techniques, we refer the reader to, e.g., \cite{DLGW,Fleischman19,DLWG,WTX}.
\end{rmk}

\begin{figure}[ht!]
\centerline{\includegraphics[trim=5.0cm 1.7cm 2.7cm 0.9cm, clip, width=14cm]{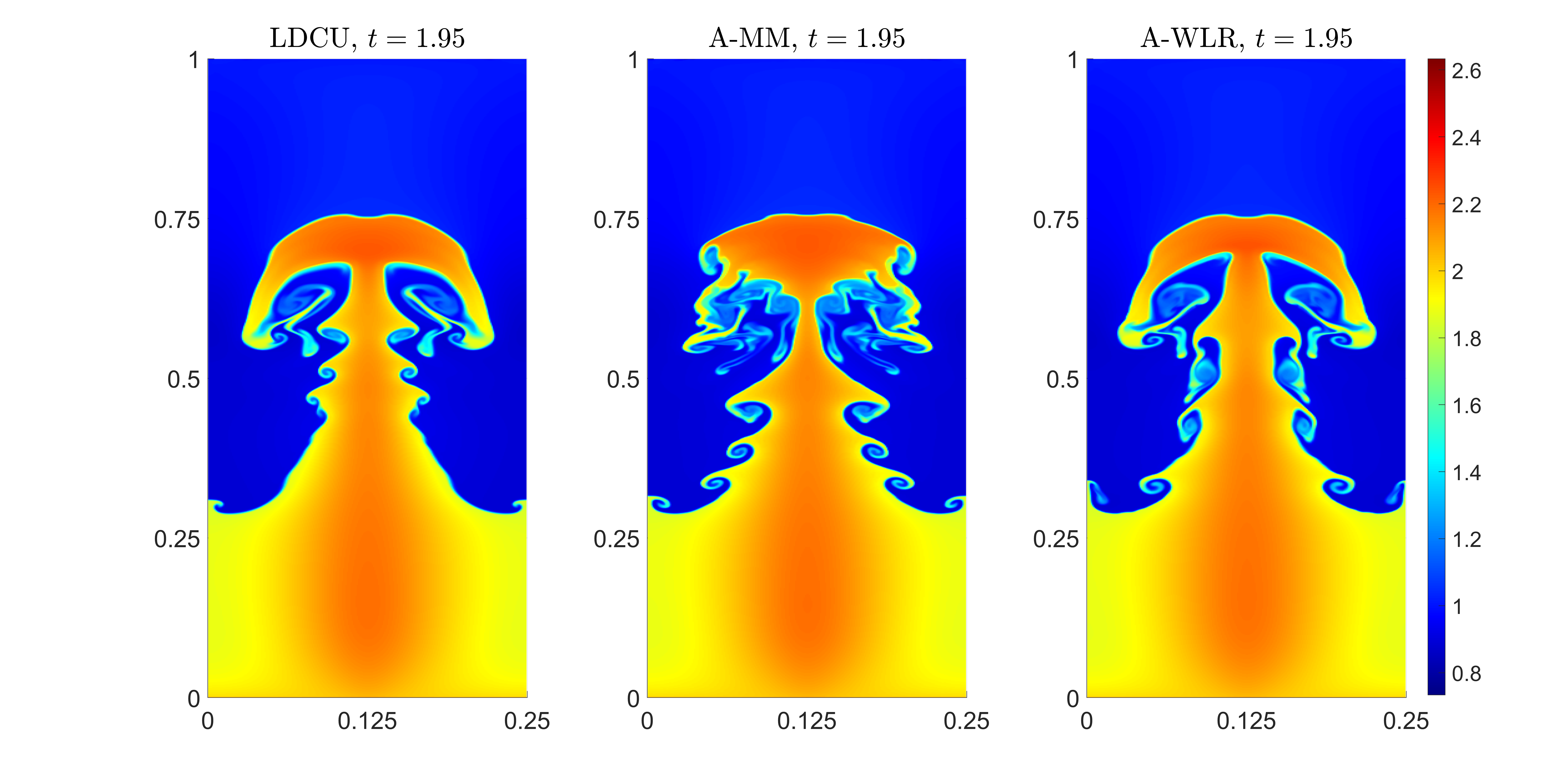}}
\vskip8pt
\centerline{\includegraphics[trim=5.0cm 1.7cm 2.7cm 0.9cm, clip, width=14cm]{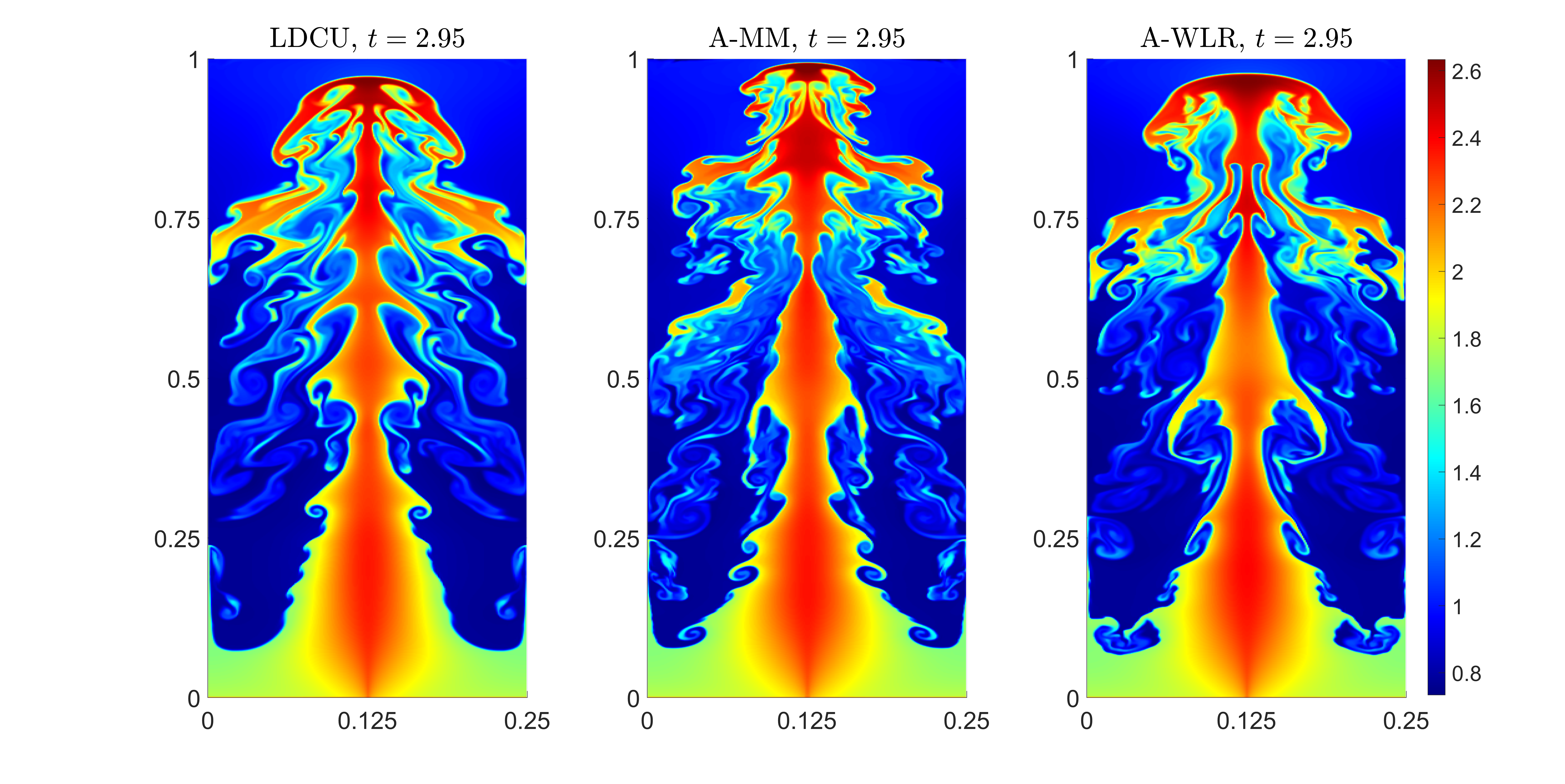}}
\caption{\sf Example 6: Density $\rho$ computed by the LDCU (left), A-MM (middle), and A-WLR (right) schemes at $t=1.95$ (top row) and 2.95
(bottom row).\label{fig10a}}
\end{figure}

\begin{figure}[ht!]
\centerline{\includegraphics[trim=4.9cm 1.7cm 4.3cm 1.0cm, clip, width=14cm]{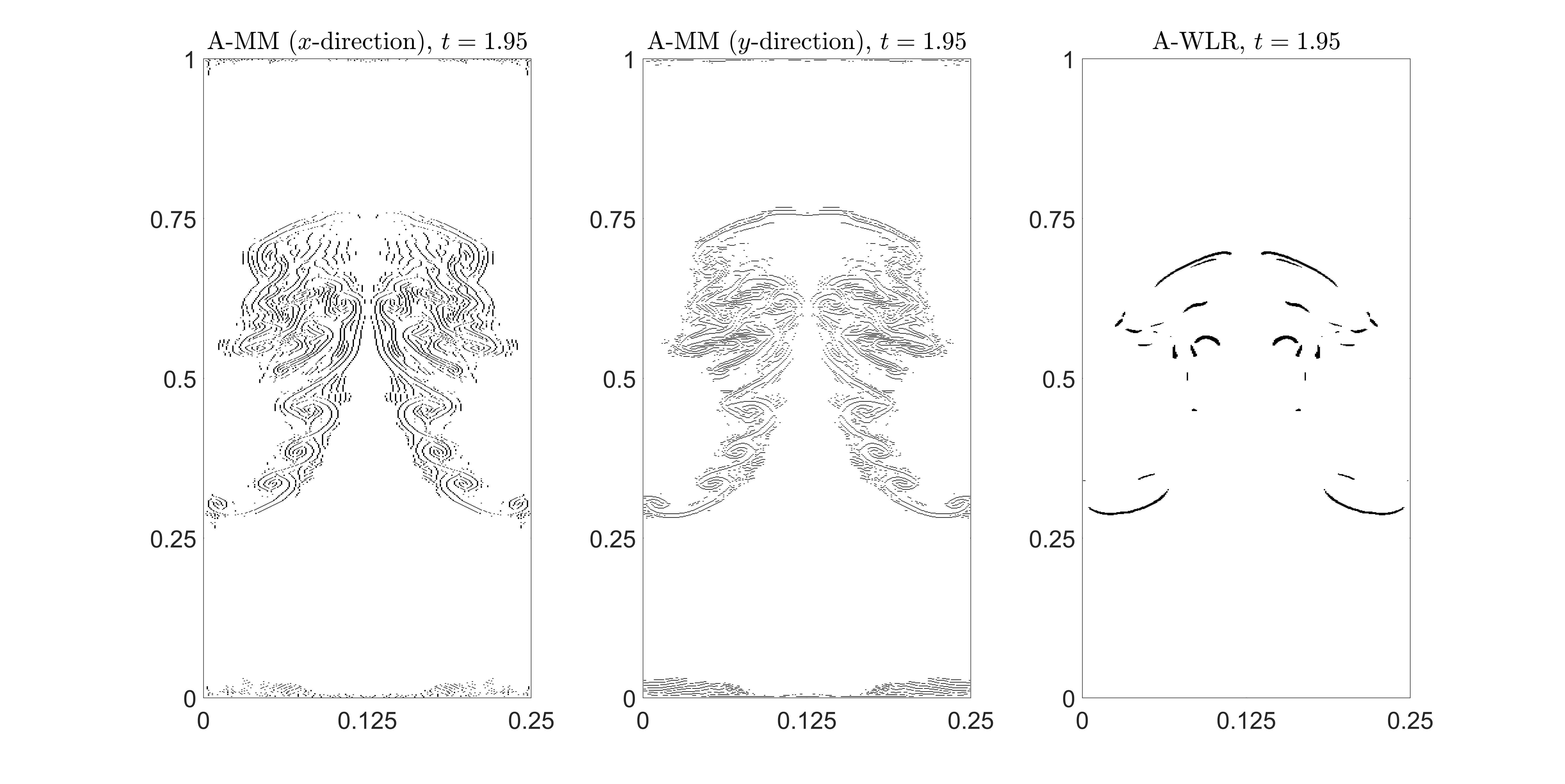}}
\vskip8pt
\centerline{\includegraphics[trim=4.9cm 1.7cm 4.3cm 1.0cm, clip, width=14cm]{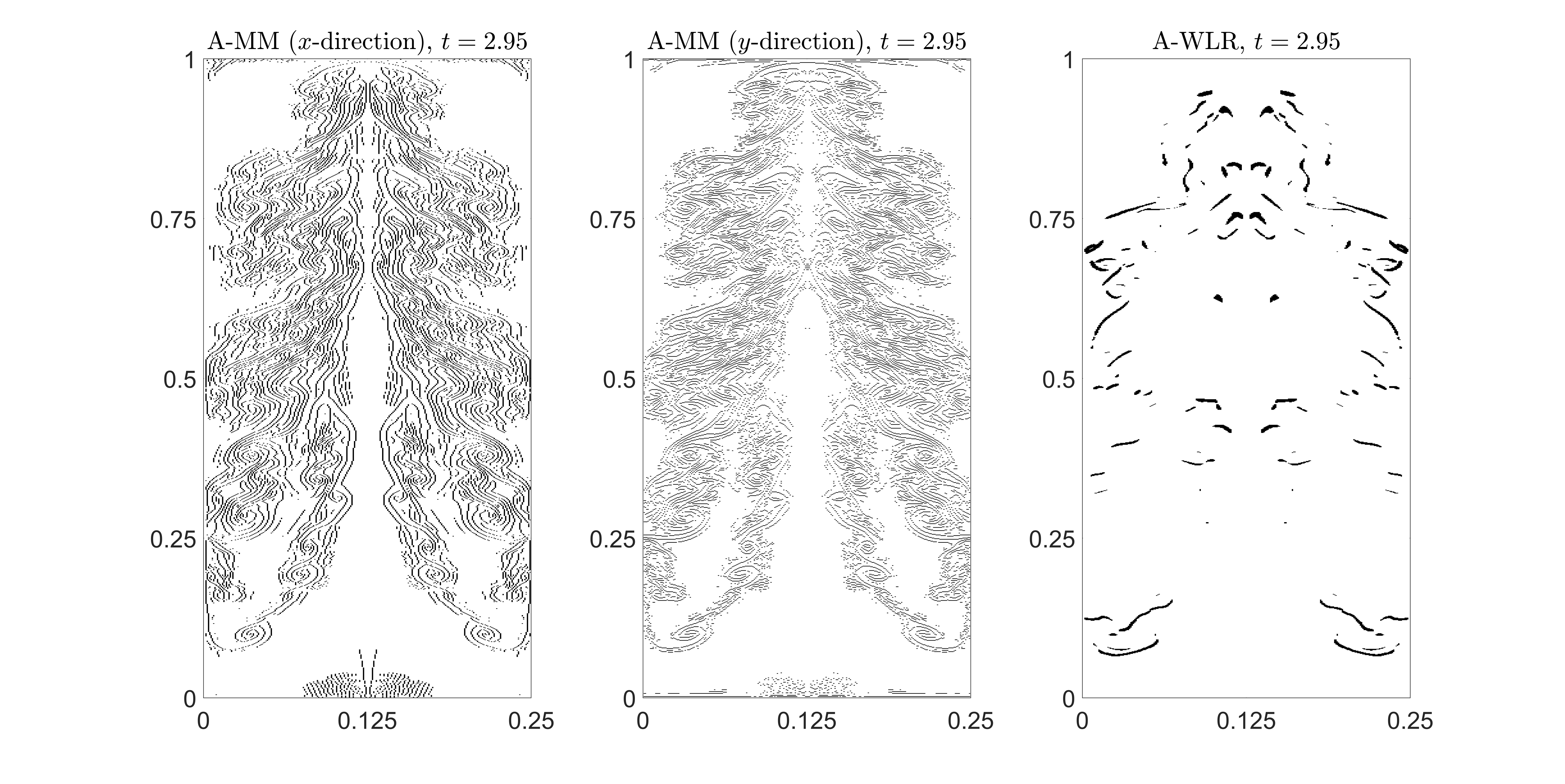}}
\caption{\sf Example 6: Areas detected as having large $x$- (left) and $y$-derivatives (middle) by the MM-based SI and the ``rough'' areas
detected by the WLR-based SI (right) at $t=1.95$ (top row) and 2.95 (bottom row).\label{fig10b}}
\end{figure}

\paragraph{Example 6---RT Instability.} In the last example taken from \cite{Shi03}, we investigate the RT instability, which is a physical
phenomenon occurring when a layer of heavier fluid is placed on top of a layer of lighter fluid. To this end, we first modify the 2-D Euler
equations of gas dynamics \eref{1.2}, \eref{3.3}--\eref{3.4} by adding the gravitational source terms acting in the positive direction of
the $y$-axis into the RHS of the system:
\begin{equation*}
\begin{aligned}
&\rho_t+(\rho u)_x+(\rho v)_y=0,\\
&(\rho u)_t+(\rho u^2 +p)_x+(\rho uv)_y=0,\\
&(\rho v)_t+(\rho uv)_x+(\rho v^2+p)_y=\rho,\\
&E_t+\left[u(E+p)\right]_x+\left[v(E+p)\right]_y=\rho v,
\end{aligned}
\end{equation*}
and then use the following initial conditions:
\begin{equation*}
(\rho,u,v,p)(x,y,0)=\begin{cases}(2,0,-0.025c\cos(8\pi x),2y+1),&y<0.5,\\(1,0,-0.025c\cos(8\pi x),y+1.5),&\mbox{otherwise},\end{cases}
\end{equation*}
where $c:=\sqrt{\gamma p/\rho}$ is the speed of sound. The solid wall boundary conditions are imposed at $x=0$ and $x=0.25$, and the
following Dirichlet boundary conditions are specified at the top and bottom boundaries:
$$
(\rho,u,v,p)(x,1,t)=(1,0,0,2.5),\quad(\rho,u,v,p)(x,0,t)=(2,0,0,1).
$$
We compute the numerical solution until the final time $t=2.95$ by the LDCU, A-MM, and A-WLR (with the adaption constant $\texttt{C}=3$)
schemes on the uniform mesh with $\dx=\dy=1/1024$ in the computational domain $[0,0.25]\times[0,1]$ and then present the numerical results
obtained at the times $t=1.95$ and 2.95 in Figure \ref{fig10a}. As one can see, the A-MM and A-WLR schemes resolve more small structures
than the LDCU scheme, which again demonstrates that using a sharper SBM limiter can produce sharper numerical results. As in the previous
examples, one can observe that the A-MM scheme uses the overcompressive SBM limiter in a larger number of cells compared with the A-WLR
scheme; see Figure \ref{fig10b}, where the indicated ``rough'' areas are plotted.
\begin{rmk}
In this example, the solution is symmetric with respect to the vertical axis $x=0.125$. In order to enforce this symmetry, we have applied
the strategy from \cite{WDGK2020}: upon completion of each time evolution step, we replace the computed cell averages $\xbar\mU_{j,k}$ with
$\widehat{\bm U}_{j,k}$, where
\begin{equation*}
\begin{aligned}
\widehat\rho_{j,k}&=\frac{\xbar\rho_{j,k}+\xbar\rho_{M-j,k}}{2},&(\widehat{\rho u})_{j,k}&
=\frac{(\xbar{\rho u})_{j,k}-(\xbar{\rho u})_{M-j,k}}{2},\\
(\widehat{\rho v})_{j,k}&=\frac{(\xbar{\rho v})_{j,k}+(\xbar{\rho v})_{M-j,k}}{2},&\widehat E_{j,k}&
=\frac{\xbar E_{j,k}+\xbar E_{M-j,k}}{2},
\end{aligned}
\end{equation*}
for all $j=1,\ldots,M$ and for all $k$. Alternative symmetry enforcement techniques can be found in, e.g.,
\cite{DLGW,DLWG,Fleischman19,WTX}.
\end{rmk}

\section{Conclusion}\label{sec5}
In this paper, we have introduced new second-order adaptive low-dissipation central-upwind schemes for the one- (1-D) and two-dimensional
(2-D) hyperbolic systems of conservation laws. The new adaptive schemes are based on the recently proposed low-dissipation central-upwind
(LDCU) fluxes and two smoothness indicators (SIs) (the minmod (MM)- and weak local residual (WLR)-based ones) used to automatically detect
``rough'' areas of the computed solutions. We then use the overcompressive SBM limiters in the ``rough'' areas and the dissipative Minmod2
limiters elsewhere to achieve higher resolution of the computed shocks and contact discontinues and, at the same time, to avoid the
staircase-like overcompressed structures in the computed results. We have applied the developed adaptive schemes to the 1-D and 2-D Euler
equations of gas dynamics and the obtained numerical results clearly demonstrate that both of the adaptive schemes outperform the original
LDCU scheme.

We have also compared the performance of the two proposed adaptive schemes. It turns out that even though the use of the WLR-based SI may be
advantageous in some examples, this SI relies on an adaption constant, which may be hard to tune: this affects the robustness of the
resulting adaptive scheme. The use of the MM-based SI, on the other hand, leads to a robust adaption strategy. Other SIs may be tested and
they may turn out to be even more robust and sharp, but we leave this study for the future work.

\section*{Acknowledgments}
The work of A. Kurganov was supported in part by NSFC grant 12171226, and by the fund of the Guangdong Provincial Key Laboratory of
Computational Science and Material Design (No. 2019B030301001).

\appendix
\section{2-D LCD-Based Piecewise Linear Reconstruction}\label{appa}
In this appendix, we describe how to reconstruct the one-sided point values $\mU^\pm_{\jph,k}$ (the point values $\mU^\pm_{j,\kph}$ can be
computed in a similar manner and we omit the details for the sake of brevity). To this end, as in the 1-D case, we first introduce the local
characteristic variables in the neighborhood of $(x,y)=(x_\jph,y_k)$:
\begin{equation*}
\bm\Gamma_{\ell,k}=R^{-1}_{\jph,k}\xbar\mU_{\ell,k},\quad\ell=j-1,\,j,\,j+1,\,j+2,
\end{equation*}
where the matrix $R_{\jph,k}$ is such that $R^{-1}_{\jph,k}\widehat A_{\jph,k}R_{\jph,k}$ is diagonal and a locally linearized Jacobian is
$\widehat A_{\jph,k}:=A\big((\,\xbar\mU_{j,k}+\xbar\mU_{j+1,k})/2\big)$.

Equipped with the values $\bm\Gamma_{j-1,k}$, $\bm\Gamma_{j,k}$, $\bm\Gamma_{j+1,k}$, and $\bm\Gamma_{j+2,k}$, we compute
\begin{equation*}
(\bm\Gamma_x)_{j,k}=\phi^{\rm SBM}_{\theta,\tau}\left(\frac{\bm\Gamma_{j+1,k}-\bm\Gamma_{j,k}}{\bm\Gamma_{j,k}-\bm\Gamma_{j-1,k}}\right)
\frac{\bm\Gamma_{j,k}-\bm\Gamma_{j-1,k}}{\dx},
\end{equation*}
and
\begin{equation*}
(\bm\Gamma_x)_{j+1,k}=\phi^{\rm SBM}_{\theta,\tau}\left(\frac{\bm\Gamma_{j+2,k}-\bm\Gamma_{j+1,k}}{\bm\Gamma_{j+1,k}-\bm\Gamma_{j,k}}\right)
\frac{\bm\Gamma_{j+1,k}-\bm\Gamma_{j,k}}{\dx},
\end{equation*}
where the SBM function, defined in \eref{2.9}, is applied in the component-wise manner. We then use these slopes to evaluate
$$
\bm\Gamma^-_{\jph,k}=\bm\Gamma_{j,k}+\frac{\dx}{2}(\bm\Gamma_x)_{j,k}\quad\mbox{and}\quad
\bm\Gamma^+_{\jph,k}=\bm\Gamma_{j+1,k}-\frac{\dx}{2}(\bm\Gamma_x)_{j+1,k},
$$
and finally obtain the corresponding point values of $\mU$ by
\begin{equation*}
\mU^\pm_\jph=R_{\jph,k}\bm\Gamma^\pm_{\jph,k}.
\end{equation*}
\begin{rmk}
The matrices $R_{\jph,k}$ and $R^{-1}_{\jph,k}$ for the 2-D Euler equation of gas dynamics \eref{3.3}--\eref{3.4} can be found in
\cite[Appendix C]{CCHKL_22}.
\end{rmk}

\bibliographystyle{siam}
\bibliography{Chu-Kurganov}
\end{document}